\documentclass[11pt,x11names]{elsarticle} 
\usepackage[a4paper,top=1in,bottom=1in,left=0.8in,right=0.8in]{geometry}
\usepackage{amsmath,mathtools}
\usepackage{amsthm,amssymb}
\theoremstyle{definition}
\newtheorem{defn}{Definition}
\allowdisplaybreaks
\usepackage{subcaption}
\usepackage{multirow}
\usepackage{tabu}
\usepackage[ruled]{algorithm2e}
\newtheorem{theorem}{Theorem}
\newtheorem{corollary}{Corollary}
\newtheorem{example}{Example}
\usepackage[x11names]{xcolor}
\usepackage{tikz-network}

\usepackage{parskip}
\usepackage{hyperref}
\hypersetup{
    colorlinks=true,
    linkcolor=blue,
    filecolor=magenta,      
    urlcolor=cyan,
    citecolor=blue,
}
\usepackage[authormarkup=none,commentmarkup=margin,xcolor=x11names]{changes}
\setcommentmarkup{\footnote{{\bf #2 \arabic{authorcommentcount}:} \ifthenelse{\isColored}{\color{authorcolor}}{}#1}}
\definechangesauthor[color=red]{AE}

\definechangesauthor[color=purple]{AK}

\usepackage{todonotes}
\long\def\remark[#1]#2{\todo[inline,author=#1]{\emph{#2}}}
\journal{Computers \& Operations Research}
\makeatletter
\def\ps@pprintTitle{%
 \let\@oddhead\@empty
 \let\@evenhead\@empty
 \def\@oddfoot{}%
 \let\@evenfoot\@oddfoot}
\makeatother

\biboptions{authoryear}

\begin{document}

\begin{frontmatter}

\title{New commodity representations for multicommodity network flow problems: An application to the fixed-charge network design problem}

\author[math]{Ahmad Kazemi}
\corref{mycorrespondingauthor}
\cortext[mycorrespondingauthor]{Corresponding author}
\ead{ahmad.kazemi@monash.edu}
\author[it]{Pierre Le Bodic}
\ead{pierre.lebodic@monash.edu}
\author[math]{Andreas T. Ernst}
\ead{andreas.ernst@monash.edu}
\author[uq]{Mohan Krishnamoorthy}
\ead{m.krishnamoorthy@uq.edu.au}

\address[math]{School of Mathematics, Monash University, Australia}
\address[it]{Faculty of Information Technology, Monash University, Australia}
\address[uq]{School of Information Technology and Electrical Engineering, The University of Queensland, Australia}

\begin{abstract}
When solving hard multicommodity network flow problems using an LP-based approach, the number of commodities is a driving factor in the speed at which the LP can be solved, as it is linear in the number of constraints and variables.
The conventional approach to improve the solve time of the LP relaxation of a Mixed Integer Programming (MIP) model that encodes such an instance is to aggregate all commodities that have the same origin or the same destination. 
However, the bound of the resulting LP relaxation can significantly worsen, which tempers the efficiency of aggregating techniques. In this paper, we introduce the concept of \emph{partial aggregation} of commodities that aggregates commodities over a subset of the network instead of the conventional aggregation over the entire underlying network.
This offers a high level of control on the trade-off between size of the aggregated MIP model and quality of its LP bound.
We apply the concept of partial aggregation to two different MIP models for the multicommodity network design problem.
Our computational study on benchmark instances confirms that the trade-off between solve time and LP bound can be controlled by the level of aggregation, and that choosing a good trade-off can allow us to solve the original large-scale problems faster than without aggregation or with full aggregation.
\end{abstract}

\begin{keyword}
aggregation; network optimization; multicommodity network flow; LP relaxation
\end{keyword}

\end{frontmatter}

\section{Introduction}\label{sec:intro}

Aggregation techniques are useful tools to tackle the challenges of solving large-scale optimization problems. These methods offer a trade-off between the mathematical model size and the level of detail included. In Linear Programming (LP) and Mixed-Integer Programming (MIP), the model size can be reduced by aggregating constraints (rows) \citep{zipkinRow}, variables (columns) \citep{zipkinCol}, or both at the same time \citep{Rogers1991}. Aggregation techniques have been successfully applied to a wide variety of problems in network flow optimization \citep{Evans1983}, mining \citep{BOLAND20091064}, environmental planning \citep{nazari-2015}, machine learning \citep{park_klabjan_2016}, and location analysis \citep{francis_lowe_rayco_tamir_2008}, to name but a few.
An aggregated MIP may have more feasible solutions than its disaggregated version, and the same holds for their LP relaxations, which leads in general to a weaker LP bound. On the other hand, size reduction translates to shorter computing times and possibly reduced degeneracy \citep{elhallaoui_metrane_soumis_desaulniers_2008}.
The trade-off between quality of the bound and the time it takes to compute it is not always beneficial to a MIP algorithm such as the branch-and-bound.
Indeed, aggregation can allow nodes to be computed faster, but the weakening of the bound can lead to a much larger branch-and-bound tree, overall requiring more compute time.
In this paper, we propose novel aggregation approaches for multicommodity network flow problems that provide much more control over this trade-off. We will show in this paper that our new aggregation approaches can significantly reduce the model size while only minimally degrading the quality of the LP relaxation, which benefits MIP algorithms to solve large instances.

Multicommodity network flow problems seek to route a set of commodities between an origin and a destination over a network. The commodities are usually interrelated by bundling constraints such as capacity constraints and/or network design decisions. Multicommodity network flow problems have been employed together with linear objective functions, piecewise-linear objective functions \citep{Croxton2007}, network design problems with economies of scale \citep{Andrews}, network loading problem \citep{AVELLA2007103}, and network design with nodes costs \citep{Belotti2007}. 
Moreover, there are applications that can be modelled as multicommodity time-space network although they may not exhibit an explicit network structure \citep{KLIEWER20061616, KAZEMI2020}. 
To show the effectiveness of the proposed aggregation approaches, we apply these to the Multicommodity Capacitated fixed-charge Network Design problem (\emph{MCND}). However, these approaches are not restricted to MCND and are applicable to any other problem that can be formulated with a multicommodity network flow subproblem in its variables and constraints. We choose MCND for three reasons. First, this problem is easy to understand, but difficult to solve. Second, MCND is well studied, has extensive applications \citep{Magnanti}, and is the basis for many other network design problems. Third, there are well-established benchmark instances for MCND in the literature \citep{CRAINIC200173}. Although, our aim is not to improve on the state-of-the-art for solving MCND, if that was an outcome of the work we carry out, it would indeed, be beneficial. What we are looking for, instead, is a generalized framework for solving problems that can be formulated with a multicommodity network flow subproblem in its variables and constraints, of which the MCND is a significant archetype.

The input of MCND is a directed graph $G = (\mathcal{N}, \mathcal{A})$ and a set of commodities $\mathcal{K}$. Each commodity $k \in \mathcal{K}$ corresponds to a flow of $d^k$ units that must be routed from an origin $o^k \in \mathcal{N}$ to a destination $s^k \in \mathcal{N}$. For each arc $(i,j) \in \mathcal{A}$, a capacity $u_{ij}$, a per-unit-of-flow cost $c_{ij}$ and a fixed cost $f_{ij}$ are defined, all non-negative. The fixed cost is incurred if an arc is used/opened for routing flow.

%
It is well-known in the MCND literature that all commodities with the same origin (with the same destination) can be aggregated as one commodity with an origin and multiple destinations (with multiple origins and a destination) \citep{Chouman2017}. This leads to two formulations with the same optimal integer solution but with different LP bound qualities: (a) a disaggregated formulation (\emph{DA}) that considers commodities based on origin-destination pairs; (b) a fully-aggregated formulation (\emph{FA})  that aggregates all commodities with the same origin into a single commodity with that origin and multiple destinations, one per aggregated commodity. The same aggregation process can be used to instead aggregate commodities that have the same destination, but, without loss of generality, we only aggregate by origins in this paper.

We employ the most common Mixed-Integer Program (MIP) model for MCND in the literature \citep{Gendron1999}. This model includes a continuous flow variable $x_{ij}^k$ that corresponds to the amount of flow of commodity $k$ on the arc $(i,j)$, and a binary variable $y_{ij}$ that determines if the arc $(i,j)$ is opened. The binary input value $o_i^k$ ($s_i^k$) is 1 if and only if the node $i\in\mathcal{N}$ is the origin (destination) of the commodity $k$, and we partition the set of neighbors of a node $i$ into the set of successor nodes $\mathcal{N}_i^+ = \{j\in\mathcal{N} \mid (i,j)\in \mathcal{A}\}$, and the set of predecessor nodes $\mathcal{N}_i^- = \{j\in\mathcal{N} \mid (j,i)\in \mathcal{A}\}$. The MCND MIP model is as follows:
\begin{subequations}
\begin{align}
  \noalign{\noindent\bf Disaggregated Formulation (DA):}
	\min \hspace{10pt} & \sum_{k \in \mathcal{K}}\sum_{(i,j) \in \mathcal{A}} c_{ij}x_{ij}^k + \sum_{(i,j)\in \mathcal{A}} f_{ij} y_{ij} \label{eq:obj}\\
    s.t: \hspace{9pt} & \sum_{j\in \mathcal{N}_i^+} x_{ij}^k - \sum_{j \in \mathcal{N}_i^-} x_{ji}^k = (o_i^k-s_i^k)d^k &&\forall\ k \in \mathcal{K},\ i \in \mathcal{N} \label{eq:flow_conserv}\\
    &\sum_{k \in \mathcal{K}} x_{ij}^k \leq u_{ij}y_{ij} &&\forall\ (i,j)\in \mathcal{A} \label{eq:capacity}\\
    &x_{ij}^k \leq d^k y_{ij} && \forall\ k \in \mathcal{K},\ (i,j) \in \mathcal{A} \label{eq:stong_ineq}\\
    &x_{ij}^k \geq 0 &&\forall\ (i,j)\in \mathcal{A} \label{eq:x_sign}\\
    &y_{ij} \in \{0,1\} &&\forall\ (i,j) \in \mathcal{A} \label{eq:y_binary}
\end{align}
\end{subequations}
Objective function~\eqref{eq:obj} minimizes the total cost, including the flow and fixed costs. Equation~\eqref{eq:flow_conserv} is the flow conservation constraint. Constraint~\eqref{eq:capacity} enforces the capacity limits of arcs. Constraint~\eqref{eq:stong_ineq} links the binary and continuous decision variables in a disaggregated form of Constraint~\eqref{eq:capacity}.
Constraint~\eqref{eq:stong_ineq} is not necessary for the validity of the MIP, however, it significantly improves the tightness of the LP relaxation bound \citep{Gendron1999}. This constraint is called a \emph{Strong Inequality} (SI) in the literature \citep{Chouman2017} and we adopt the same name throughout. The decision variables and their domains are defined in Constraints~\eqref{eq:x_sign}-\eqref{eq:y_binary}. This model corresponds to the DA formulation. The FA formulation, which is based on the full aggregation approach that aggregates all commodities with a same origin, is similarly developed by modifications in the commodity definition and summing up the flow conservation constraints~\eqref{eq:flow_conserv} and the strong inequalities~\eqref{eq:stong_ineq} over the aggregated commodities. 

\cite{Gendron1999} showed that, in practice, the LP relaxation of DA provides tight bounds for MCND. On the other hand, FA is a smaller formulation with significantly shorter LP relaxation computing time, but looser bounds. The DA formulation has received more attention particularly to develop MIP-based heuristics \citep{crainic_gendron_hernu_2004, RODRIGUEZMARTIN2010575, GENDRON201870}. The FA formulation has been preferred in cutting plane algorithms by \cite{bienstock_1995,Bienstock1996}, \cite{bienstock_1998} and \cite{Christian} because of the shorter computing time of its LP relaxation. \cite{RARDIN199395} studied the relationship between the aggregated and the disaggregated commodity representation for the uncapacitated fixed-charge network design with one origin and multiple destinations. The aggregated version of this problem includes one commodity while its disaggregated version includes one commodity per destination with nonzero demand. The authors derive an improved formulation by adding a family of  so-called ``dicut'' inequalities to the FA formulation that has the same LP relaxation polyhedral as the DA formulation for this problem. These inequalities are specific to the problem they considered. \cite{Chouman2017} are the first who evaluated the effect of commodity representation, either disaggregated or fully-aggregated, on the performance of a Branch-and-Bound algorithm (B\&B) with a specialized cutting plane algorithm. They showed that in general, the DA formulation outperforms the FA formulation over the instances with many commodities, but the FA formulation is more effective for instances with few commodities. 

The model selection in the literature follows an ``all-or-nothing'' approach that employs either the DA or the FA formulation based on the application and the algorithm considered. This approach sacrifices one aspect, LP relaxation computational difficulty or LP relaxation bound quality, for the other aspect.  In this paper, we propose a spectrum of formulations that lie in within the range of the DA and FA formulations in this trade-off. This allows to select a formulation with the desired level of compromise between the model computational difficulty and the LP bound tightness. These formulations are based on \emph{partial} aggregations of commodities instead of a disaggregated or a conventional fully aggregated approach.

We make four main contributions in this paper. First, we introduce new commodity representations for the multicommodity network flow problems that make the partial aggregations possible. Second, we propose a base partially-aggregated formulation. We improve the base model by two types of tightening constraints and propose a heuristic to construct effective partial aggregations that gives a proper trade-off between the LP bound quality and computing time. Third, we study the polyhedra of the LP relaxations of the proposed formulations and their relation with polyhedra of the LP relaxation of the FA and DA formulations. Fourth, we empirically investigate the LP relaxations of the proposed MIPs and compare them with the FA and DA formulations. We also evaluate the performance of the MIP algorithms over the formulations on benchmark instances.

The paper is organized as follows. Section~\ref{sec:commodity} extends the conventional aggregation schemes to define \emph{dispersions} of commodities and partial aggregations as generalization. Moreover, a MIP is designed for these new types of aggregations. Section~\ref{sec:formulations} proposes two types of specialized tightening constraints for the partially-aggregated formulation, and presents a heuristic algorithm to construct effective partial aggregations. These formulations result in a spectrum of formulations within the range of FA and DA formulations in terms of the LP relaxation tightness and computing time. Section~\ref{sec:polyhed} studies the polyhedra of the proposed formulations and compares them with the FA and the DA formulations. Section~\ref{sec:comp} presents an extensive computational study of the proposed formulations and their LP relaxations. Section~\ref{sec:conclusion} concludes the paper.

\section{Commodity definitions and aggregation levels}\label{sec:commodity}

Recall that in the input, each commodity $k \in \mathcal{K}$ is defined by three attributes; an origin $o^k$, a destination $s^k$, and a demand of $d^k$. 
The representation of commodities in the DA formulation is congruous with the definition of commodities in the input. 
In particular, there is a one-to-one mapping between commodities of the input and commodities in the DA formulation. 
However, this need not be the case, as illustrated in the FA formulation, in which a ``\emph{commodity}'' is an aggregation of the input commodities that share the same origin.

\begin{figure}[htb!]
    \centering
    \begin{subfigure}[t]{0.49\textwidth}
    	\centering
                \begin{tikzpicture}[multilayer=3d,scale=0.8]
\SetLayerDistance{2.1}
\Plane[x=0,y=0,width=4.5,height=2.4,layer=1,color=white,opacity=1]
\Vertex[x=.5,y=1.2,size=0.25,layer=1,color=red,shape = rectangle]{1n1}
\Vertex[x=1.5,y=1.6,size=0.25,layer=1,color=white]{1n2}
\Vertex[x=1.5,y=.8,size=0.25,layer=1,color=white]{1n3}
\Vertex[x=2.5,y=2,size=0.25,layer=1,color=white]{1n4}
\Vertex[x=2.5,y=1.2,size=0.25,layer=1,color=white]{1n5}
\Vertex[x=2.5,y=.4,size=0.25,layer=1,color=white]{1n6}
\Vertex[x=3.5,y=1.6,size=0.25,layer=1,color=red]{1n7}
\Vertex[x=3.5,y=.8,size=0.25,layer=1,color=white]{1n8}

\draw[ ->, >=stealth,color=red](1n1) -- (1n2) ;
\draw[ ->, >=stealth,color=red](1n1) -- (1n3) ;
\draw[ ->, >=stealth,color=red](1n2) -- (1n4) ;
\draw[ ->, >=stealth,color=red](1n2) -- (1n5) ;
\draw[ ->, >=stealth,color=red](1n3) -- (1n5) ;
\draw[ ->, >=stealth,color=red](1n3) -- (1n6) ;
\draw[ ->, >=stealth,color=red](1n5) -- (1n7) ;
\draw[ ->, >=stealth,color=red](1n5) -- (1n8) ;
\draw[ ->, >=stealth,color=red](1n7) -- (1n4) ;
\draw[ ->, >=stealth,color=red](1n8) -- (1n6) ;
\Text[x=4.2,y=.3,layer=1,color=red]{$k_1$}
\Plane[x=0,y=0,width=4.5,height=2.4,layer=2,color=white,opacity=1]
\Vertex[x=.5,y=1.2,size=0.25,layer=2,color=teal,shape = rectangle]{2n1}
\Vertex[x=1.5,y=1.6,size=0.25,layer=2,color=white]{2n2}
\Vertex[x=1.5,y=.8,size=0.25,layer=2,color=white]{2n3}
\Vertex[x=2.5,y=2,size=0.25,layer=2,color=white]{2n4}
\Vertex[x=2.5,y=1.2,size=0.25,layer=2,color=white]{2n5}
\Vertex[x=2.5,y=.4,size=0.25,layer=2,color=white]{2n6}
\Vertex[x=3.5,y=1.6,size=0.25,layer=2,color=white]{2n7}
\Vertex[x=3.5,y=.8,size=0.25,layer=2,color=teal]{2n8}

\draw[ ->, >=stealth,color=teal](2n1) -- (2n2) ;
\draw[ ->, >=stealth,color=teal](2n1) -- (2n3) ;
\draw[ ->, >=stealth,color=teal](2n2) -- (2n4) ;
\draw[ ->, >=stealth,color=teal](2n2) -- (2n5) ;
\draw[ ->, >=stealth,color=teal](2n3) -- (2n5) ;
\draw[ ->, >=stealth,color=teal](2n3) -- (2n6) ;
\draw[ ->, >=stealth,color=teal](2n5) -- (2n7) ;
\draw[ ->, >=stealth,color=teal](2n5) -- (2n8) ;
\draw[ ->, >=stealth,color=teal](2n7) -- (2n4) ;
\draw[ ->, >=stealth,color=teal](2n8) -- (2n6) ;
\Text[x=4.2,y=.3,layer=2,color=teal]{$k_2$}

\Plane[x=0,y=0,width=4.5,height=2.4,layer=3,color=white,opacity=1]
\Vertex[x=.5,y=1.2,size=0.25,layer=3,color=blue,shape = rectangle]{3n1}
\Vertex[x=1.5,y=1.6,size=0.25,layer=3,color=white]{3n2}
\Vertex[x=1.5,y=.8,size=0.25,layer=3,color=white]{3n3}
\Vertex[x=2.5,y=2,size=0.25,layer=3,color=blue]{3n4}
\Vertex[x=2.5,y=1.2,size=0.25,layer=3,color=white]{3n5}
\Vertex[x=2.5,y=.4,size=0.25,layer=3,color=white]{3n6}
\Vertex[x=3.5,y=1.6,size=0.25,layer=3,color=white]{3n7}
\Vertex[x=3.5,y=.8,size=0.25,layer=3,color=white]{3n8}

\draw[ ->, >=stealth,color=blue](3n1) -- (3n2) ;
\draw[ ->, >=stealth,color=blue](3n1) -- (3n3) ;
\draw[ ->, >=stealth,color=blue](3n2) -- (3n4) ;
\draw[ ->, >=stealth,color=blue](3n2) -- (3n5) ;
\draw[ ->, >=stealth,color=blue](3n3) -- (3n5) ;
\draw[ ->, >=stealth,color=blue](3n3) -- (3n6) ;
\draw[ ->, >=stealth,color=blue](3n5) -- (3n7) ;
\draw[ ->, >=stealth,color=blue](3n5) -- (3n8) ;
\draw[ ->, >=stealth,color=blue](3n7) -- (3n4) ;
\draw[ ->, >=stealth,color=blue](3n8) -- (3n6) ;
\Text[x=4.2,y=.3,layer=3,color=blue]{$k_3$}
\Plane[x=0,y=0,width=4.5,height=2.4,layer=4,color=white,opacity=1]
\Vertex[x=.5,y=1.2,size=0.25,layer=4,color=purple,shape = rectangle]{4n1}
\Vertex[x=1.5,y=1.6,size=0.25,layer=4,color=white]{4n2}
\Vertex[x=1.5,y=.8,size=0.25,layer=4,color=white]{4n3}
\Vertex[x=2.5,y=2,size=0.25,layer=4,color=white]{4n4}
\Vertex[x=2.5,y=1.2,size=0.25,layer=4,color=white]{4n5}
\Vertex[x=2.5,y=.4,size=0.25,layer=4,color=purple]{4n6}
\Vertex[x=3.5,y=1.6,size=0.25,layer=4,color=white]{4n7}
\Vertex[x=3.5,y=.8,size=0.25,layer=4,color=white]{4n8}

\draw[ ->, >=stealth,color=purple](4n1) -- (4n2) ;
\draw[ ->, >=stealth,color=purple](4n1) -- (4n3) ;
\draw[ ->, >=stealth,color=purple](4n2) -- (4n4) ;
\draw[ ->, >=stealth,color=purple](4n2) -- (4n5) ;
\draw[ ->, >=stealth,color=purple](4n3) -- (4n5) ;
\draw[ ->, >=stealth,color=purple](4n3) -- (4n6) ;
\draw[ ->, >=stealth,color=purple](4n5) -- (4n7) ;
\draw[ ->, >=stealth,color=purple](4n5) -- (4n8) ;
\draw[ ->, >=stealth,color=purple](4n7) -- (4n4) ;
\draw[ ->, >=stealth,color=purple](4n8) -- (4n6) ;
\Text[x=4.2,y=.3,layer=4,color=purple]{$k_4$}

\end{tikzpicture}
        \caption{Implied network by the DA formulation}
        \label{fig:da-net}
    \end{subfigure}
    \begin{subfigure}[t]{0.49\textwidth}
    	\centering
        \begin{tikzpicture}[multilayer=3d,scale=0.8]
\SetLayerDistance{2.7}
\Plane[x=0,y=0,width=4.5,height=2.6,layer=0,color=white,opacity=1,NoBorder]

\Plane[x=0,y=0,width=4.5,height=2.8,layer=1,color=white,opacity=1]
\Vertex[x=.5,y=1.6,size=0.25,layer=1,color=black,shape = rectangle]{1n1}
\Vertex[x=1.5,y=2,size=0.25,layer=1,color=white]{1n2}
\Vertex[x=1.5,y=1.2,size=0.25,layer=1,color=white]{1n3}
\Vertex[x=2.5,y=2.4,size=0.25,layer=1,color=blue]{1n4}
\Vertex[x=2.5,y=1.6,size=0.25,layer=1,color=white]{1n5}
\Vertex[x=2.5,y=.8,size=0.25,layer=1,color=purple]{1n6}
\Vertex[x=3.5,y=2,size=0.25,layer=1,color=red]{1n7}
\Vertex[x=3.5,y=1.2,size=0.25,layer=1,color=teal]{1n8}

\draw[ ->, >=stealth,color=black](1n1) -- (1n2) ;
\draw[ ->, >=stealth,color=black](1n1) -- (1n3) ;
\draw[ ->, >=stealth,color=black](1n2) -- (1n4) ;
\draw[ ->, >=stealth,color=black](1n2) -- (1n5) ;
\draw[ ->, >=stealth,color=black](1n3) -- (1n5) ;
\draw[ ->, >=stealth,color=black](1n3) -- (1n6) ;
\draw[ ->, >=stealth,color=black](1n5) -- (1n7) ;
\draw[ ->, >=stealth,color=black](1n5) -- (1n8) ;
\draw[ ->, >=stealth,color=black](1n7) -- (1n4) ;
\draw[ ->, >=stealth,color=black](1n8) -- (1n6) ;
\Text[x=3.2,y=.3,layer=1,color=black]{$\{k_1,k_2,k_3,k_4\}$}

\end{tikzpicture}        \caption{Implied network by the FA formulation}
        \label{fig:fa-net}
    \end{subfigure}
    \caption{Network representation of the DA and FA formulations for an instance with 4 commodities, all originating from the square node in the figure, on a graph with 8 nodes and 10 arcs.}
    \label{fig:network-agg}
\end{figure}
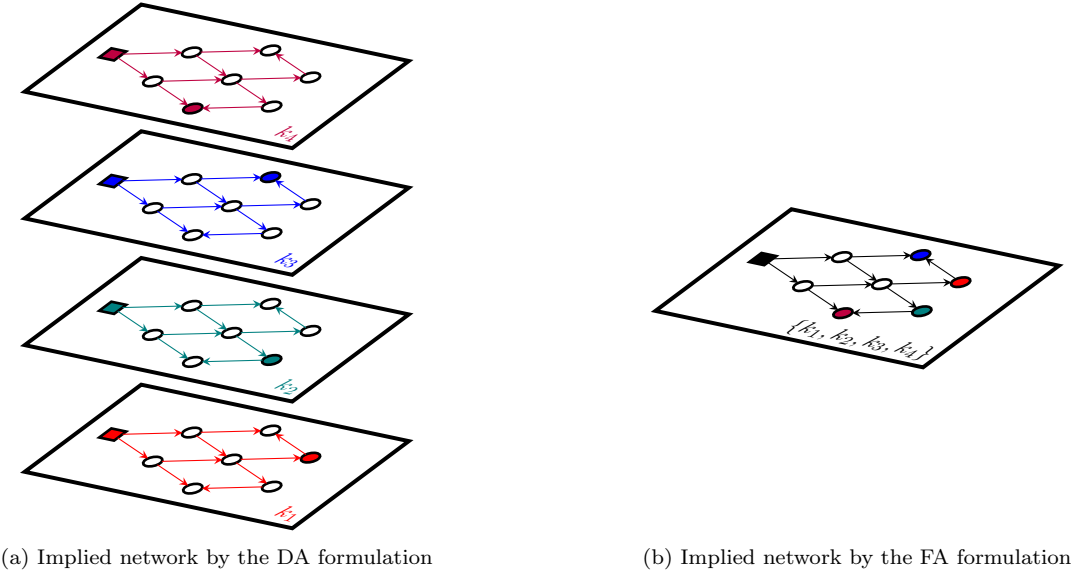

Alternative commodity definitions actually correspond to different network structures that represent the problem. This results in various sizes of the MIP model that encodes the problem since the nodes and arcs of the network correspond to constraints and variables of the MIP model. These MIPs have the same optimal integer solution but different LP relaxations.  In the DA formulation, the flow conservation constraints~\eqref{eq:flow_conserv}  describe the routes of each commodity on the graph $G$ independently from other commodities. This means that the DA formulation essentially considers a hypothetical network \emph{layer} for each commodity. This hypothetical layer is a copy of the original graph $G$ on which the demand of its corresponding commodity flow.
On the other hand, in the FA formulation, all layers related to the aggregated commodities are merged into a single layer on which the total demand of the aggregated commodities flow. Figure~\ref{fig:network-agg} shows the network representations of the DA and FA formulations for an instance with 4 commodities. In this instance, the original graph $G$ consists of 8 nodes and 10 arcs, and all commodities are originated from the square node. Figure~\ref{fig:da-net} shows the network representation of the DA formulation. In this figure, the opaque nodes represent the origin and destination of the corresponding commodity of each layer. As there are 4 commodities, network representation of the DA formulation includes 4 layers. These layers are merged into one layer by the FA formulation as shown in Figure~\ref{fig:fa-net}. This layer includes an origin (square) node for the aggregated commodities and four opaque (circle) nodes as destinations of the aggregated commodities. 

In this section, we introduce new commodity definitions that give a control on the aggregation level that could range from the disaggregated approach to the fully aggregated approach.  This enables us to develop a spectrum of formulations that fill in the gap between the DA and the FA formulations in terms of the trade-off between the computational difficulty and bound quality of the LP relaxation of a formulation. The intention of our approach is to allow groups of commodities to start in an aggregated form at a shared origin and to be disaggregated as they move towards their individual destination. Therefore, the set of commodities included in a same group varies over the network. For this purpose we define a \emph{dispersion} of commodities.
A dispersion is defined as a set of commodities and specifies how each of its commodities is aggregated on each arc.

  \begin{defn}[dispersion] A dispersion $b$ of commodities on a directed graph $G=(\mathcal{N},\mathcal{A})$ is defined by the following information:
    \begin{enumerate}
    \item A set $\mathcal{K}_b\subseteq \mathcal{K}$ of commodities that share the common origin $o_b\in\mathcal{N}$.
    \item The set of destination nodes $\mathcal{S}_b\subseteq \mathcal{N}$ ($\mathcal{S}_b=\{s^k: k\in \mathcal{K}_b\}$).
    \item On each arc $(i,j) \in \mathcal{A}$, a partition of $\mathcal{K}_b$ into a subset $\mathcal{K}_b^{ij}$ of commodities aggregated on that arc and its complement $\mathcal{D}_{b}^{ij}$ of commodities which are disaggregated on that arc.
    \end{enumerate}
    \label{def:cluster}
  \end{defn}

\begin{figure}[htb!]
    \centering
        \begin{tikzpicture}
        \node[rectangle,draw,thick,fill=black] at (-.5,1.6) (n1) {};
        \node[circle,draw,thick] at (1.5,2.5) (n2) {};
        \node[circle,draw,thick] at (1.5,.7) (n3) {};
        \node[circle,draw,thick,fill=blue] at (3.5,3.4) (n4) {};
        \node[circle,draw,thick] at (3.5,1.6) (n5) {};
        \node[circle,draw,thick,fill=purple] at (3.5,-.2) (n6) {};
        \node[circle,draw,thick,fill=red] at (5.5,2.5) (n7) {};
        \node[circle,draw,thick,fill=teal] at (5.5,.7) (n8) {};
        \path[ ->, >=stealth, thick]   (n1)  edge[] (n2);
        \path[ ->, >=stealth, thick]   (n1)  edge[] (n3);
        \path[ ->, >=stealth, thick]   (n2)  edge[bend left=20] (n4);
        \path[ ->, >=stealth, thick]   (n2)  edge[bend right=20,blue] (n4);
        \path[ ->, >=stealth, thick]   (n2)  edge[bend left=20] (n5);
        \path[ ->, >=stealth, thick]   (n2)  edge[bend right=20,red] (n5);
        \path[ ->, >=stealth, thick]   (n3)  edge[bend left=20] (n5);
        \path[ ->, >=stealth, thick]   (n3)  edge[bend right=20,teal] (n5);
        \path[ ->, >=stealth, thick]   (n3)  edge[bend left=20] (n6);
        \path[ ->, >=stealth, thick]   (n3)  edge[bend right=20,purple] (n6);
        \path[ ->, >=stealth, thick]   (n5)  edge[bend left=20] (n7);
        \path[ ->, >=stealth, thick]   (n5)  edge[bend right=20,red] (n7);
        \path[ ->, >=stealth, thick]   (n5)  edge[bend left=20] (n8);
        \path[ ->, >=stealth, thick]   (n5)  edge[bend right=20,teal] (n8);
        \path[ ->, >=stealth, thick]   (n7)  edge[] (n4);
        \path[ ->, >=stealth, thick]   (n8)  edge[] (n6);
        \draw[very thick] (-1,-1) rectangle (6, 3.7);
        \node[] at (4.8,-.7) {\small $\{k_1,k_2,k_3,k_4\}$};
		\end{tikzpicture}
    \caption{An example dispersion based on the aggregated commodities of Figure~\ref{fig:fa-net}. In this figure, colored arcs represent the disaggregation of the corresponding commodity from the group on that arc.} 
    \label{fig:dispersion}
\end{figure}
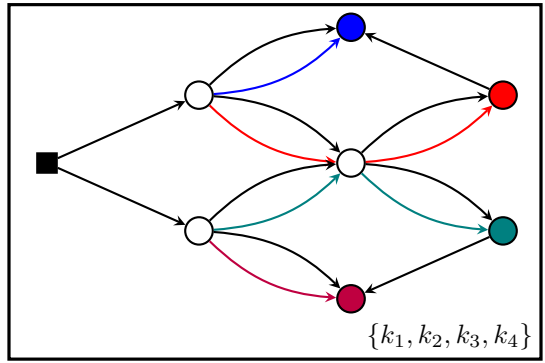

A dispersion is essentially similar to a fully-aggregated commodity but with disaggregation of some commodities on some arcs.
Suppose a commodity is disaggregated from the group on the arc $(i,j)$.
Such disaggregation is represented by adding an extra arc with the same origin and destination as the arc $(i,j)$ but exclusive to the flow of the disaggregated commodity.
Therefore, the network representation of a dispersion is similar to the network representation of a fully-aggregated commodity but with some extra arcs dedicated to the flow of the separated/disaggregated commodities, and there is one network layer for a dispersion. Figure~\ref{fig:dispersion} shows an example dispersion based on the fully-aggregated commodity shown in Figure~\ref{fig:fa-net}. In this figure, extra arcs that represent the flow of the disaggregated commodities are shown in colors. We will discuss the selection of such arcs for the commodities in Section~\ref{sec:k-path}. Black arcs correspond to the flow of the commodities that are not disaggregated from the group on that arc.

To represent an instance of the problem by dispersions, we have to define a set of dispersions in which there is at least one dispersion that includes each commodity. 
This ensures that all possible flow assignments of each commodity to all arcs of the network are covered.
We next define partial aggregations in which there is a unique dispersion for each commodity.

  \begin{defn}[partial aggregation]
    A partial aggregation for the MCND with $G=(\mathcal{N},\mathcal{A})$, original commodities $\mathcal{K}$ and demand $d^k$, consists of a set $\mathcal{B}$ of dispersions, where for every commodity $k\in\mathcal{K}$ there exists a unique $b \in \mathcal{B}$ such that $k\in \mathcal{K}_b$.
    \label{def:par-agg}
  \end{defn}


As defined in Definitions~\ref{def:cluster} and \ref{def:par-agg}, partial aggregations must satisfy three conditions: (a) aggregated commodities must always share a common origin, (b) each commodity belongs to a single dispersion, and (c) a commodity is either disaggregated on an arc or it is grouped with others of the group (that are also not disaggregated on this arc). An alternative definition is to consider dispersions based on same destinations or a mixed approach in which some dispersions include commodities with a same origin and others include commodities with a same destination. It should be noted that these definitions could of course be extended in a natural way to allow different types of dispersions and aggregations.

Conventional disaggregation and full aggregation are special cases of partial aggregations:
  \begin{description}
  \item[DA] The original, fully disaggregated version (as a special partial aggregation case) is defined by $\mathcal{B}$, where $|\mathcal{B}|=|\mathcal{K}|$, each dispersion $b \in \mathcal{B}$ includes only one commodity $\mathcal{K}_b = \{k_b\}$, and $\cup_{b \in \mathcal{B}} \mathcal{K}_b = \mathcal{K} $. Therefore, $\mathcal{S}_b = \{i\in \mathcal{N}\mid s^{k_b}_i=1\},\ |\mathcal{S}_b|=1$ for all $b\in \mathcal{B}$, and $\mathcal{K}_b^{ij}= \emptyset \And \mathcal{D}_b^{ij} = \{k_b\} : \forall b \in \mathcal{B},\  (i,j) \in \mathcal{A}$.
  \item[FA] Let $\Tilde{\mathcal{N}}$ be the set of nodes from which at least one commodity originates. The fully aggregated version is defined by $\mathcal{B}$, where $|\mathcal{B}|=|\mathcal{\Tilde{N}}|$ and  $\cup_{b\in \mathcal{B}} \{o_b\} = \Tilde{\mathcal{N}}$. Therefore,
    $\mathcal{S}_b=\{i \in \mathcal{N}\mid \exists k\in\mathcal{K} : s^k_i=1 , o^k=o_b\}$ for all $b \in \mathcal{B}$,  and $\mathcal{K}_b^{ij}=\mathcal{S}_b \And \mathcal{D}_b^{ij} = \emptyset : \forall b \in \mathcal{B},\ {(i,j)} \in \mathcal{A}$.
  \end{description}
  
The network implied by DA has one complete replica of graph $G$ for each commodity $k\in\mathcal{K}$, which results in a network with $\vert \mathcal{K} \vert$ layers. The network implied by FA has one complete replica of graph $G$ for each dispersion $n\in \mathcal{\Tilde{N}}$, which results in a network with $\vert \mathcal{\Tilde{N}} \vert$ layers.
In both cases, each layer is an exact copy of the graph $G$. However, the network layer that represents a dispersion includes some extra arcs that correspond to the disaggregation of commodities on those arcs. In particular, such a layer includes: (a) a complete copy of the nodes of graph $G$, (b) one copy of arc $(i,j) \in \mathcal{A}$ for which $\mathcal{K}_b^{ij} \neq \emptyset$, and (c) one copy of arc $(i,j) \in \mathcal{A}$ for each $k \in \mathcal{D}_b^{ij}$ to represent disaggregated arcs. An example network layer corresponding to a dispersion is shown in Figure~\ref{fig:dispersion}. Here we define a network layer that represent a dispersion as \emph{dispersion layer}.

\begin{defn}[dispersion layer]\label{def:dispersion-layer}
  A dispersion layer that represents a dispersion $b$ includes a set of nodes $\mathcal{C}_b \coloneqq \mathcal{N}$ and a set of arcs $\mathcal{A}_b^{ij}$ for each $(i,j) \in \mathcal{A}$. Each arc $a \in \mathcal{A}_b^{ij}$ corresponds to a commodity set of the set $\mathcal{G}_b^{ij} = \mathcal{K}_b^{ij} \cup \{ \{k\}\mid k \in \mathcal{D}_b^{ij}\}$ by a one-to-one mapping. Each set $D \in \mathcal{G}_b^{ij}$ implies that in the dispersion layer corresponding to dispersion $b$ there is a copy of arc $(i,j)$ that corresponds to the commodity set $D$. 
\end{defn}

A partial aggregation consists of a set of dispersions. Similarly, the network implied by a partial aggregation, called \emph{partial aggregation network}, is a set of dispersion layers, where each layer corresponds to a dispersion $b \in \mathcal{B}$.  
  A partial aggregation network that corresponds to to the partial aggregation $\mathcal{B}$ includes $|\mathcal{B}|$ layers. A partial aggregation network has $|\mathcal{B}||\mathcal{N}|$ nodes and $\sum_{b \in \mathcal{B}} \sum_{(i,j)\in \mathcal{A}} |\mathcal{G}_b^{ij}|$ arcs. The DA and FA networks have respectively $|\mathcal{K}||\mathcal{N}|$ and $|\mathcal{\Tilde{N}}||\mathcal{N}|$ nodes and $|\mathcal{K}||\mathcal{A}|$ and $|\mathcal{\Tilde{N}}||\mathcal{A}|$ arcs.

To have a formulation incorporating a partial aggregation and its corresponding partial aggregation network, we can now modify the DA formulation \eqref{eq:obj}--\eqref{eq:y_binary} by
  \begin{enumerate}
  \item Summing (\ref{eq:flow_conserv}) over $k\in \mathcal{K}_b$ for all $b\in \mathcal{B},\ i \in \mathcal{C}_b$, resulting in $\vert \mathcal{B} \vert \vert \mathcal{N} \vert$ constraints.
  \item Summing (\ref{eq:stong_ineq}) over $k\in D$ for all $(i,j)\in \mathcal{A}, b \in \mathcal{B}, D\in \mathcal{G}_b^{ij}$, resulting in between $\vert \mathcal{A} \vert \vert \mathcal{B} \vert$ and $\vert \mathcal{A} \vert \vert \mathcal{K} \vert$ constraints.
  \item Replacing any sum over variables $x_{ij}^k$ in the constraints of the resulting model by a new variable $x^D_{ij}=\sum_{k\in D}x^k_{ij}$ that corresponds to aggregated flow of a set of commodities $D \in \mathcal{G}_b^{ij}$ for some 
  $(i,j)\in \mathcal{A}$ and $b \in \mathcal{B}$.
  \end{enumerate}
This gives a new formulation that in general has both fewer constraints and fewer variables than DA.
  \begin{subequations}
\begin{align}
  \noalign{\noindent\bf Partially-Aggregated Formulation (PA):}
  \min \hspace{10pt} & \sum_{(i,j) \in \mathcal{A}} c_{ij}\sum_{b \in \mathcal{B}}\sum_{D\in \mathcal{G}_b^{ij}}
                       x_{ij}^D + \sum_{(i,j)\in \mathcal{A}} f_{ij} y_{ij} \label{eq:ag:obj}\\
  s.t: \hspace{9pt} & \sum_{j\in \mathcal{N}_i^+}\sum_{D\in\mathcal{G}_b^{ij}} x_{ij}^D -
                      \sum_{j \in \mathcal{N}_i^-}\sum_{D\in\mathcal{G}_b^{ji}} x_{ji}^D = \sum_{k\in \mathcal{K}_b} (o_i^k-s_i^k)d^k &&\forall b \in \mathcal{B},\ i \in \mathcal{C}_b \label{eq:ag:flow_conserv}\\
    &\sum_{b\in \mathcal{B}}\sum_{D \in \mathcal{G}_b^{ij}} x_{ij}^D \leq u_{ij}y_{ij} &&\forall\ (i,j)\in \mathcal{A} \label{eq:ag:capacity}\\
    & x_{ij}^D \leq \left(\sum_{k\in D}d^k\right) y_{ij} && \forall\ (i,j) \in \mathcal{A},\ b \in \mathcal{B},\ D \in\mathcal{G}_b^{ij} \label{eq:ag:SI}\\
    &x_{ij}^D \geq 0 &&\forall\ (i,j)\in \mathcal{A},\ b\in\mathcal{B},\ D\in\mathcal{G}_b^{ij}  \label{eq:ag:x_sign}\\
    &y_{ij} \in \{0,1\} &&\forall\ (i,j) \in \mathcal{A} \label{eq:ag:y_binary}
\end{align}
In the PA formulation, equation~\eqref{eq:ag:flow_conserv} conserves the aggregated flow of commodities in a same dispersion at each node of the network, whether they are disaggregated from the group or not. 
Note that if $\sum_{k\in D}d^k >u_{ij}$ for $D\in\mathcal{G}_b^{ij}$, as normally happens in the FA formulation, then the SI constraints \eqref{eq:ag:SI} become redundant.

\section{Improving the partially-aggregated formulation}\label{sec:formulations}

While (partial) aggregation makes the MIP model smaller, it weakens the formulation in terms of the LP relaxation bound. To combat this we employ two strategies that are discussed in this section, namely (a) adding tightening constraints and (b) making judicious choices about which commodities to aggregate/disaggregate for each arc. If the coefficient of variable $y_{ij}$ in the corresponding SI \eqref{eq:ag:SI} is large enough, the structure of a dispersion layer allows a commodity to flow on arc $(i,j)$ on the aggregated arc even if this commodity has been disaggregated on arc $(i,j)$. The tightening constraints we introduce in this section prevent the commodities from flowing on such arcs. We employ two approaches for this. The first approach in Section~\ref{sec:ineq-mip} adds a type of inequalities that are similar to flow conservation constraints. The second approach in Section~\ref{sec:eq-mip} modifies the partial aggregation network by adding some artificial nodes and arcs that translates to more constraints and variables.
However, the extra constraints can be expressed as equality flow conservation constraints, which both tighten the bound and, in practice, speed-up the computation of the LP. In Section~\ref{sec:k-path}, we propose a heuristic that uses a K-shortest paths algorithm to construct dispersions of the commodities and form a partial aggregation that provides a good trade-off between model size and LP bound. Here we are particularly interested in dispersions where the disaggregated arcs for each commodity form a connected network leading to the destination node of that commodity.

\subsection{Partially-aggregated formulation with inequality tightening constraints} \label{sec:ineq-mip}
We start with an example that shows a fractional solution of the LP relaxation of the PA formulation, and propose valid inequalities that can cut it. The example is shown for the incoming and outgoing flow of the node $i$ in Figure~\ref{fig:ineq fig a}. This node is part of the larger network of the problem, but the other parts are not shown as this node is enough to illustrate the example.  This node has two incoming arcs from nodes 1 and 2, and two outgoing arcs to nodes 3 and 4. Consider the dispersion $b$ of commodities $\mathcal{K}_b = \{k_1, k_2, k_3\}$. Origin and destination nodes of commodities of set $\mathcal{K}_b$ are not shown in this figure. However, in our examples we consider feasible solutions that flow of some commodities of this set pass the node $i$ as an intermediate node.
In this figure, we assume that all commodities of $\mathcal{K}_b$ are grouped together on these arcs, therefore $\mathcal{G}_b^{1i}=\mathcal{G}_b^{2i}=\mathcal{G}_b^{i3}=\mathcal{G}_b^{i4}=\mathcal{K}_b$ and there is no disaggregated arc.
As mentioned in Section~\ref{sec:commodity}, each arc of the network corresponds to a commodity set. In Figure~\ref{fig:ineq fig}, the corresponding commodity set to each arc is labeled on it. In Figure~\ref{fig:ineq fig a}, commodities are grouped on all arcs shown, and hence, all arcs correspond to the commodity set $\mathcal{K}_b$. Now, suppose that we disaggregate the commodity $k_1$ on arcs $(1,i)$ and $(i,3)$. Such disaggregation requires the network to have two additional arcs, shown in blue in Figure~\ref{fig:ineq fig b}.

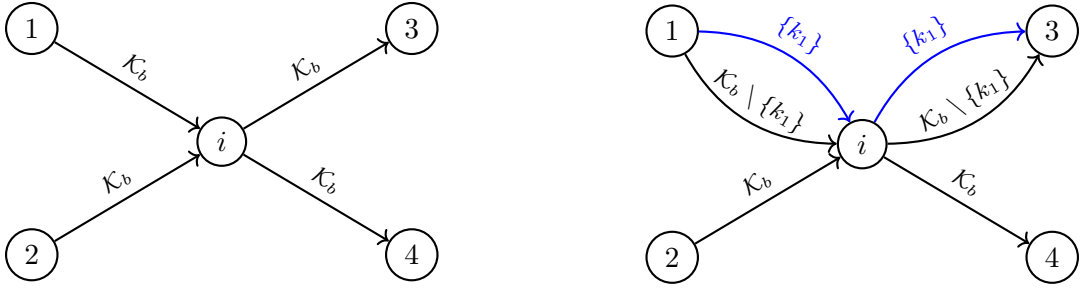
\begin{figure}[htb!]
    \centering
    \begin{subfigure}[h]{0.49\textwidth}
    	\centering
        \begin{tikzpicture}
        \node[circle,draw,thick] at (0,0) (i) {$i$};
        \node[circle,draw,thick] at (-2.5,1.5) (1) {$1$};
        \node[circle,draw,thick] at (-2.5,-1.5) (2) {$2$};
        \node[circle,draw,thick] at (2.5,1.5) (3) {$3$};
        \node[circle,draw,thick] at (2.5,-1.5) (4) {$4$};
        \path[->, thick]   (1)  edge  node[above,sloped]{\footnotesize$\mathcal{K}_b$} (i);
        \path[->, thick]   (2)  edge  node[above,sloped]{\footnotesize$\mathcal{K}_b$}(i);
        \path[->, thick]   (i)  edge node[above,sloped]{\footnotesize$\mathcal{K}_b$} (3);
        \path[->, thick]   (i)  edge node[above,sloped]{\footnotesize$\mathcal{K}_b$} (4);
		\end{tikzpicture}
        \caption{Node $i$ in case of aggregating all commodities $\mathcal{K}_b$ on the arcs related to this node}
        \label{fig:ineq fig a}
    \end{subfigure}
    \begin{subfigure}[h]{0.49\textwidth}
    	\centering
        \begin{tikzpicture}
        \node[circle,draw,thick] at (0,0) (i) {$i$};
        \node[circle,draw,thick] at (-2.5,1.5) (1) {$1$};
        \node[circle,draw,thick] at (-2.5,-1.5) (2) {$2$};
        \node[circle,draw,thick] at (2.5,1.5) (3) {$3$};
        \node[circle,draw,thick] at (2.5,-1.5) (4) {$4$};
        \path[->, thick]   (1)  edge[bend right=30] node[above,sloped]{\footnotesize$\mathcal{K}_b \setminus \{k_1\}$} (i);
        \path[->, blue,thick]   (1)  edge[bend left=30] node[above,sloped]{\footnotesize$\{k_1\}$} (i);
        \path[->, thick]   (2)  edge  node[above,sloped]{\footnotesize$\mathcal{K}_b$}(i);
        \path[->, thick]   (i)  edge[bend right=30] node[above,sloped]{\footnotesize$\mathcal{K}_b\setminus \{k_1\}$} (3);
        \path[->, blue,thick]   (i)  edge[bend left=30] node[above,sloped]{\footnotesize$\{k_1\}$} (3);
        \path[->, thick]   (i)  edge node[above,sloped]{\footnotesize$\mathcal{K}_b$} (4);
		\end{tikzpicture}
        \caption{Node $i$ in case of disaggregating the commodity $k_1$ on arcs $(1,i)$ and $(i,3)$}
        \label{fig:ineq fig b}
    \end{subfigure}
    \caption{An example of fully aggregated and partial aggregation network. The example shows node $i$ in the dispersion layer $b$ with two incoming and two outgoing arcs. Each arc is labeled with its corresponding commodity set.}
    \label{fig:ineq fig}
\end{figure}

As the flow conservation constraint~\eqref{eq:ag:flow_conserv} considers all commodity sets $\mathcal{G}_b^{ij}$ related to a dispersion on an arc in one equation. Therefore, the flow of the disaggregated commodity $k$ may flow on the arc that is related to the commodity set $\mathcal{K}_b^{ij}$ even if $k \notin \mathcal{K}_b^{ij}$, in case the coefficient for variable $y_{ij}$ in the corresponding SI constraint \eqref{eq:ag:SI} of $\mathcal{K}_b^{ij}$ is large enough, which degrades the LP bound.  For instance, in Figure~\ref{fig:ineq fig b} the flow of commodity $k_1$ can arrive to node $i$ by the blue arc $(1,i)$ and leave it by the black arc $(i,3)$. However, the related commodity set of the black arc $(i,3)$ does not include commodity $k_1$. Let $d^{k_1} = 1$. Consider Example~\ref{exmp:n1} which is a feasible flow for the PA formulation over node $i$. 

\begin{example}\label{exmp:n1}
This is an example feasible solution for the PA formulation for the flow on incoming and outgoing arcs of node $i$ in Figure~\ref{fig:ineq fig b}: $x_{1i}^{\{k_1\}} =x_{i3}^{\mathcal{K}_b\setminus\{k_1\}} =  1$ and flow on other incoming and outgoing arcs of node $i$ equals to zero.
\end{example}

Although Example~\ref{exmp:n1} is a feasible solution for the PA, we want to prevent such flows as this flow arrives at node $i$ by an arc dedicated to $k_1$. However, it leaves node $i$ by an arc related to a commodity set that does not include $k_1$.  To prevent such flow assignments we add two inequalities~\eqref{eq:path-ineq-1} and~\eqref{eq:path-ineq-2}.
\begin{align}
	& \sum_{j \in \mathcal{N}_i^+}\sum_{D \in \mathcal{G}_{b}^{ij}:k\in D} x_{ij}^D - \sum_{j \in \mathcal{N}_i^-}\sum_{D \in \mathcal{G}_{b}^{ji}:D= \{k\}} x_{ji}^D \geq (o_i^k-s_i^k)d^k &&\forall\ b \in \mathcal{B},\ k\in\mathcal{K}_b,\ i \in \mathcal{C}_{b} \label{eq:path-ineq-1}\\ 
	& \sum_{j \in \mathcal{N}_i^+}\sum_{D \in \mathcal{G}_{b}^{ij}:D=\{k\}} x_{ij}^D - \sum_{j \in \mathcal{N}_i^-}\sum_{D \in \mathcal{G}_{b}^{ji}:k\in D} x_{ji}^D \leq (o_i^k-s_i^k)d^k &&\forall\ b \in \mathcal{B},\ k\in\mathcal{K}_b,\ i \in \mathcal{C}_{b} \label{eq:path-ineq-2}
\end{align}
  \end{subequations}
Intuitively, Constraint~\eqref{eq:path-ineq-1}, called  \emph{forward labeling inequality}, states that any flow that arrives at node $i$ labelled as commodity $k$ (i.e. with a $x_{ij}^D$ variable where $D=\{k\}$) must leave either as commodity $k$ or as an aggregated commodity that includes $k$. Example~\ref{exmp:n1} violates this constraint, as in this case, the left-hand side of this of constraint is $-1$ while its right-hand side is zero. Constraint~\eqref{eq:path-ineq-2}, called  \emph{backward labeling inequality}, creates the symmetric requirement on arriving arcs for an individual commodity $k$ leaving node $i$. A feasible solution of the PA formulation for this instance as $x_{1i}^{\mathcal{K}_b\setminus\{k_1\}} = x_{i3}^{\{k_1\}} =  1$ and the flow on other incoming and outgoing arcs as zero violates Constraint~\eqref{eq:path-ineq-2}. In the rest of the paper, we consider \eqref{eq:ag:obj}-\eqref{eq:path-ineq-2} as the \emph{partially aggregated formulation with forward and backward labeling inequalities}, and refer to it as \textbf{PAi}.

\subsection{Partially-aggregated formulation with equality tightening constraints} \label{sec:eq-mip}

In this section, we present another approach to tackle the issue raised in Section~\ref{sec:ineq-mip} for the PA formulation. In this approach, instead of adding inequalities, we modify the partial aggregation network to differentiate the incoming flow to a node labeled as a disaggregated commodity from the aggregated incoming flow. Consider the node $i$ of the dispersion layer $b$ with some incoming and outgoing arcs (aggregated and disaggregated). We add some artificial nodes and arcs in this node to be able to differentiate the flow of disaggregated commodities. As this approach adds extra nodes and arcs to the partial aggregation network, its corresponding MIP would be larger than the MIP of PAi for a same partial aggregation $\mathcal{B}$. However, its merit is that the added constraints are equality flow conservation constraints which keeps the multicommodity network structure of the problem.

As formulated, there is one flow conservation constraint per dispersion layer and per node.
This means that a node of a dispersion layer does not differentiate the incoming and outgoing flow by the disaggregated arcs from the incoming and outgoing flow by the aggregated arcs, which implies that the flow of one commodity can be routed on arcs that do not represent that commodity.
Here we propose to decompose every node into a small gadget that will track which commodity can arrive and leave on which arc, and restrict the possibilities of incorrect routing.

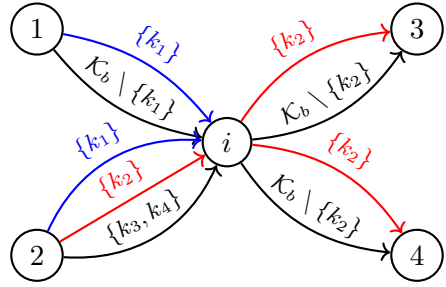
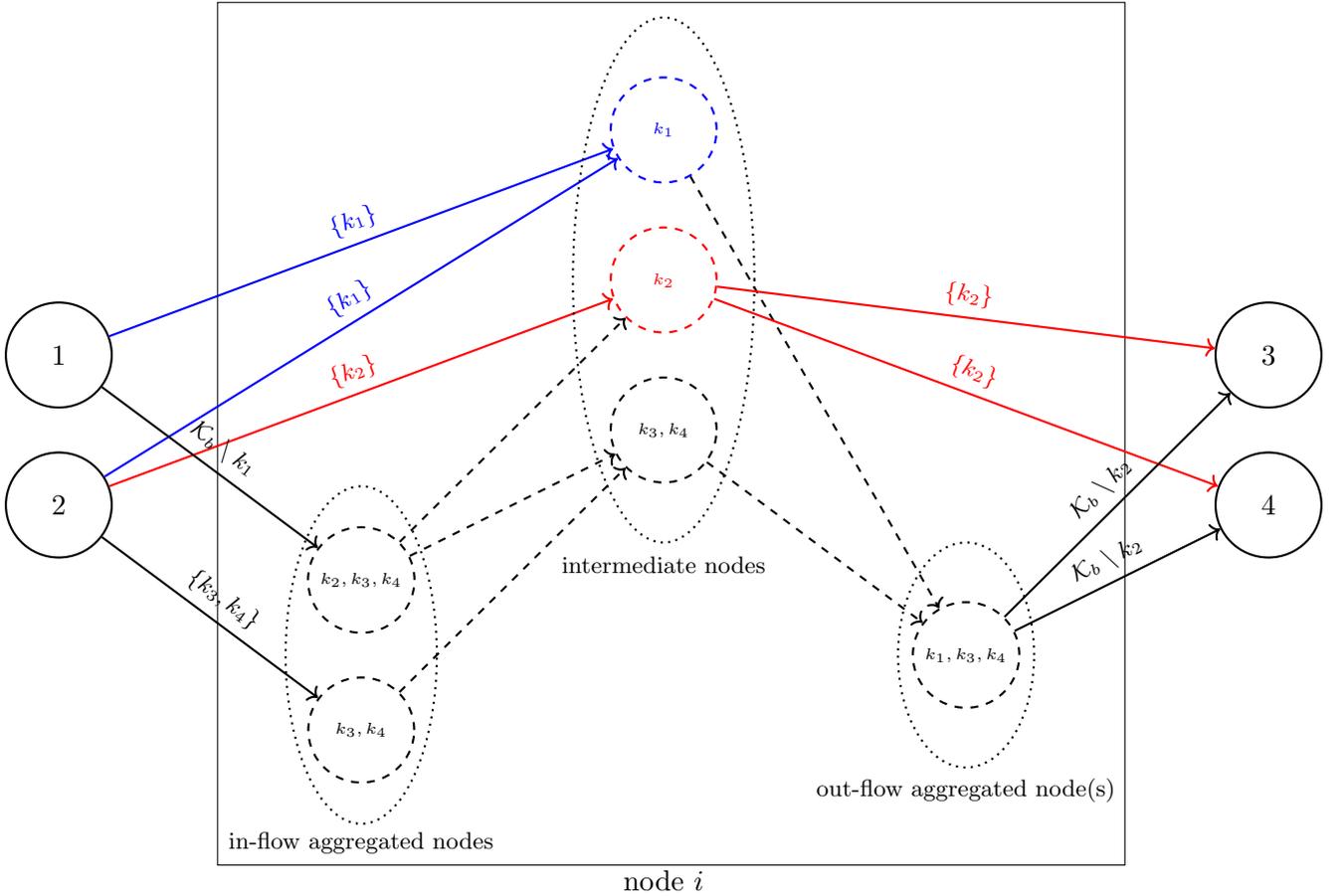
\begin{figure}[htb!]
    \centering
    \begin{subfigure}[h]{0.99\textwidth}
    	\centering
        \begin{tikzpicture}
        \node[circle,draw,thick] at (0,0) (i) {$i$};
        \node[circle,draw,thick] at (-2.5,1.5) (1) {$1$};
        \node[circle,draw,thick] at (-2.5,-1.5) (2) {$2$};
        \node[circle,draw,thick] at (2.5,1.5) (3) {$3$};
        \node[circle,draw,thick] at (2.5,-1.5) (4) {$4$};
        \path[->, thick,blue]   (1)  edge[bend left=20]  node[above,sloped]{\footnotesize$\{k_1\}$} (i);
        \path[->, thick]   (1)  edge[bend right=20]  node[above,sloped]{\footnotesize$\mathcal{K}_b\setminus\{k_1\}$} (i);
        \path[->, thick,blue]   (2)  edge[bend left=35]   node[above,sloped]{\footnotesize$\{k_1\}$}(i);
        \path[->, thick,red]   (2)  edge  node[above,sloped]{\footnotesize$\{k_2\}$}(i);
        \path[->, thick]   (2)  edge[bend right=35]   node[above,sloped]{\footnotesize$\{k_3,k_4\}$}(i);
        \path[->, thick,red]   (i)  edge[bend left=25]  node[above,sloped]{\footnotesize$\{k_2\}$} (3);
        \path[->, thick]   (i)  edge[bend right=25]  node[above,sloped]{\footnotesize$\mathcal{K}_b\setminus\{k_2\}$} (3);
        \path[->, thick,red]   (i)  edge[bend left=25]  node[above,sloped]{\footnotesize$\{k_2\}$} (4);
        \path[->, thick]   (i)  edge[bend right=25]  node[above,sloped]{\footnotesize$\mathcal{K}_b\setminus\{k_2\}$} (4);
		\end{tikzpicture}
        \caption{Node $i$ of dispersion layer $b$. Node $i$ in the original graph $G$ has two incoming arcs and two outgoing arcs. Dispersion $b$ includes commodities $\mathcal{K}_b=\{k_1,k_2,k_3,k_4\}$. Commodity $k_1$ is disaggregated on arcs $(1,i)$ and $(2,i)$. Commodity $k_2$ is disaggregated on arcs $(2,i)$, $(i,3)$, and $(i,4)$. The corresponding commodity set of each arc is labeled on it.}
        \label{fig:eq fig a}
    \end{subfigure}\\
    \begin{subfigure}[h]{0.99\textwidth}
    	\centering
        \begin{tikzpicture}
        \draw (-5.9,-9.8) rectangle (6.1, 1.7);
        \node[] at (0,1) () {};
        \node[circle,draw,thick,dashed,blue,minimum size = 1.4cm] at (0,0) (m-1) {\tiny $k_1$};
        \node[circle,draw,thick,dashed,red,minimum size = 1.4cm] at (0,-2) (m-2) {\tiny $k_2$};
        \node[circle,draw,thick,dashed,black,minimum size = 1.4cm] at (0,-4) (m-3) {\tiny $k_3,k_4$};
        \node[circle,draw,thick,dashed,black,minimum size = 1.4cm] at (-4,-6) (in-1) {\tiny $k_2,k_3,k_4$};
        \node[circle,draw,thick,dashed,black,minimum size = 1.4cm] at (-4,-8) (in-2) {\tiny $k_3,k_4$};
        \node[circle,draw,thick,dashed,black,minimum size = 1.4cm] at (4,-7) (out-1) {\tiny $k_1,k_3,k_4$};
        \path[->, thick,dashed]   (in-1)  edge (m-2);
        \path[->, thick,dashed]   (in-1)  edge (m-3);
        \path[->, thick,dashed]   (in-2)  edge (m-3);
        \path[->, thick,dashed]   (m-1)  edge (out-1);
            \path[->, thick,dashed]   (m-3)  edge (out-1);
        \node[circle,draw,thick,minimum size = 1.4cm] at (-8,-3) (n1) {$1$};
        \node[circle,draw,thick,minimum size = 1.4cm] at (-8,-5) (n2) {$2$};
        \node[circle,draw,thick,minimum size = 1.4cm] at (8,-3) (n3) {$3$};
        \node[circle,draw,thick,minimum size = 1.4cm] at (8,-5) (n4) {$4$};
        \path[->, thick,blue]   (n1)  edge  node[above,sloped]{\footnotesize$\{k_1\}$}(m-1);
        \path[->, thick,blue]   (n2)  edge  node[above,sloped]{\footnotesize$\{k_1\}$}(m-1);
        \path[->, thick,red]   (n2)  edge  node[above,sloped]{\footnotesize$\{k_2\}$}(m-2);
        \path[->, thick,black]   (n1)  edge  node[above,sloped]{\footnotesize$\mathcal{K}_b\setminus k_1$}(in-1);
        \path[->, thick,black]   (n2)  edge  node[above,sloped]{\footnotesize$\{k_3,k_4\}$}(in-2);
        \path[->, thick,red]   (m-2)  edge  node[above,sloped]{\footnotesize$\{k_2\}$}(n3);
        \path[->, thick,red]   (m-2)  edge  node[above,sloped]{\footnotesize$\{k_2\}$}(n4);
        \path[->, thick]   (out-1)  edge  node[above,sloped]{\footnotesize$\mathcal{K}_b\setminus k_2$}(n3);
        \path[->, thick]   (out-1)  edge  node[above,sloped]{\footnotesize$\mathcal{K}_b\setminus k_2$}(n4);
        \node[shape=ellipse,thick,draw,dotted, style={minimum height=4.5cm,minimum width=2cm}] at (-4,-7) (e1) {};
        \node[]  at (-4, -9.5) {\footnotesize in-flow aggregated nodes};
        \node[shape=ellipse,thick,draw,dotted, style={minimum height=7cm,minimum width=2.4cm}] at (0,-2) (e2)  {};
        \node[]  at (0, -5.8) {\footnotesize intermediate nodes};
        \node[shape=ellipse,thick,draw,dotted, style={minimum height=3cm,minimum width=1.8cm}] at (4,-7) (e3) {};
        \node[]  at (4, -8.8) {\footnotesize out-flow aggregated node(s)};
        \node[]  at (0, -10) {node $i$};
		\end{tikzpicture}
        \caption{Modified node $i$ as the larger rectangular node with added artificial nodes and arcs shown by dashed circles and arrows.}
        \label{fig:eq fig b}
    \end{subfigure}
    \caption{An illustrative example for network modification for the partially-aggregated formulation with equality tightening constraints}
    \label{fig:eq fig}
\end{figure}

As an example, consider the node $i$ in the dispersion layer $b$ shown in Figure~\ref{fig:eq fig a}. Dispersion $b$ includes 4 commodities as $\mathcal{K}_b=\{k_1,k_2,k_3,k_4\}$. Node $i$ is part of a larger network, but other nodes are not shown as this node is enough to show the examples we intend. Moreover, origin and destination nodes of commodities of set $\mathcal{K}_b$ are not shown in this figure. Node $i$ in the original graph $G$ has two incoming arcs $(1,i)$ and $(2,i)$ and two outgoing arcs $(i,3)$ and $(i,4)$. Commodity $k_1$ is disaggregated from the group on arcs $(1,i)$ and $(2,i)$. Commodity $k_2$ is disaggregated from the group on arcs $(2,i)$, $(i,3)$, and $(i,4)$. The original node $i$ in the partial aggregation network is shown in Figure~\ref{fig:eq fig a}. One possible incorrect flow route for this example is a flow arriving at node $i$ by the red arc between nodes $2$ and $i$ which corresponds to the commodity $k_2$, but leaving the node by the black arc between nodes $i$ and $3$ which does not include commodity $k_2$ in its corresponding set.

The gadget we propose to tackle the raised issue considers the node $i$ as a larger node that itself includes a set of some artificial nodes, where each of such nodes correspond to a specific commodity set. 
We define the group of artificial nodes for each node and their connection as in fact a network modification that adds artificial nodes and arcs inside a node of a dispersion layer. 
For each node $i$ and dispersion layer $b$, we create three groups of artificial nodes:
\begin{itemize}
    \item in-flow aggregated nodes: one node for each distinct commodity set that incoming aggregated arcs of node $i$ represent. For node $i$ in Figure~\ref{fig:eq fig a}, we add two in-flow aggregated nodes, one corresponds to the commodity set $\mathcal{K}_b\setminus\{k_1\}$ and the other one corresponds to the commodity set $\{k_3,k_4\}$.
    \item intermediate nodes: a node for each commodity $k \in \mathcal{K}_b$ that has at least one incident disaggregated arc to node $i$ and one node for the rest of commodities of dispersion $b$. For node $i$ in Figure~\ref{fig:eq fig a}, we add three intermediate nodes. Two nodes are added because two commodities, $k_1 \text{ and } k_2$, have incident disaggregated arcs to node $i$. The third node is added to represent the rest of commodities of set $\mathcal{K}_b$, i.e. $\{k_3,k_4\}$.
    \item out-flow aggregated nodes: one node for each distinct commodity set that outgoing aggregated arcs of node $i$ represent. For node $i$ in Figure~\ref{fig:eq fig a}, we add one out-flow aggregated node as both outgoing aggregated arcs correspond to the same commodity set, i.e. $\mathcal{K}_b\setminus k_2$.
\end{itemize}
Then, there is an artificial arc from any in-flow aggregated node to any intermediate node if and only if they represent at least one common commodity.
Similarly, there is an artificial arc from any intermediate node to any out-flow aggregated node if and only if they represent at least one common commodity.
Such artificial arcs facilitate the modelling as they describe the relation of artificial nodes added.
They have no flow cost or opening cost, hence their flow are not present in the objective function and there is no associated binary variable as they are always open. Original arcs of the partial aggregation network remain same. However, they must be connected to the relevant artificial node. Disaggregated arcs connect to the intermediate node of the corresponding disaggregated commodity. Incoming (outgoing) aggregated arcs connect to the corresponding in-flow (out-flow) aggregated node. 
 Node $i$ in Figure~\ref{fig:eq fig a} is modified and illustrated as the larger rectangular node in Figure~\ref{fig:eq fig b}, where artificial nodes and arcs are shown by dashed circles and arrows. Moreover, each group of artificial nodes inside the node $i$ is shown in a dotted ellipse for clarity.

Specifically, consider the following sets for each $b\in \mathcal{B}$ and $i \in \mathcal{C}_b$.
\begin{center}
	\begin{tabular}{cp{16cm}}
		$\mathcal{L}_b^i$ &Set of commodities of dispersion $b$ that have at least one incident disaggregated arc to node $i$. For node $i$ in Figure~\ref{fig:eq fig}, we have $\mathcal{L}_b^i = \{k_1,k_2\}$. \\
		$\mathcal{M}_b^i$ &Set of commodity sets that each added intermediate node for node $i$ of dispersion $b$ represents. For node $i$ in Figure~\ref{fig:eq fig}, we have $\mathcal{M}_b^i = \{\{k_1\},\{k_2\},\{k_3,k_4\}\}$.\\
		$\check{\mathcal{T}}_b^i$ &Set of distinct commodity sets that each incoming aggregated arc of node $i$ corresponds to. For node $i$ of Figure~\ref{fig:eq fig}, we have $\check{\mathcal{T}}_b^i = \{\{k_2,k_3,k_4\}, \{k_3,k_4\}\}$.\\
		$\hat{\mathcal{T}}_b^i$ &Set of distinct commodity sets that each outgoing aggregated arc of node $i$ corresponds to. For node $i$ of Figure~\ref{fig:eq fig}, we have $\hat{\mathcal{T}}_b^i = \{\{k_1,k_3,k_4\}\}$.
	\end{tabular}
\end{center}
\setcounter{table}{0}

Formally, these sets for the node $i$ of the dispersion layer $b$ are obtained by equations~\eqref{eq:pae-dis-nodes}-\eqref{eq:pae-out-agg-nodes}. Artificial nodes and arcs are added inside the node $i$ of the dispersion layer $b$ only if $\mathcal{L}_b^i \neq \emptyset$.
\begin{subequations}
\begin{align}
&\mathcal{L}_b^i = \{k\in \mathcal{K}_b \mid \exists\ (j,i) \in \mathcal{N}_i^-:k\in \mathcal{D}_b^{ji}\ or\ \exists\ (i.j) \in \mathcal{N}_i^+:k\in \mathcal{D}_b^{ij} \} \label{eq:pae-dis-nodes}\\ 
&\mathcal{M}_b^i = \{\{k\}\mid k \in \mathcal{L}_b^i\} \cup \{\mathcal{K}_b\setminus\mathcal{L}_b^i\} \label{eq:pae-m-nodes}\\
&\nonumber \check{\mathcal{T}}_b^i = \{C_b^{ik}\subseteq\mathcal{K}_b\mid k=1,\dots,|\check{\mathcal{T}}_b^i|, C_b^{ik}\neq \emptyset\}\\ &\hspace{25pt}\text{such that: } \cup_{C\in \check{\mathcal{T}}_b^i} C = \cup_{j \in \mathcal{N}_i^-}\{\mathcal{K}_b^{ji}\} \text{ and } C\cap D = \emptyset \;  \forall\ C\neq D \in \check{\mathcal{T}}_b^i \label{eq:pae-in-agg-nodes}\\
&\nonumber \hat{\mathcal{T}}_b^i = \{C_b^{ik}\subseteq\mathcal{K}_b\mid k=1,\dots,|\hat{\mathcal{T}}_b^i|, C_b^{ik}\neq \emptyset\}\\ &\hspace{25pt}\text{such that: } \cup_{C\in \hat{\mathcal{T}}_b^i} C = \cup_{j \in \mathcal{N}_i^+}\{\mathcal{K}_b^{ij}\} \text{ and } C\cap D = \emptyset \;  \forall\ C\neq D \in \hat{\mathcal{T}}_b^i \label{eq:pae-out-agg-nodes}
\end{align}
\end{subequations}


Consider Example~\ref{exmp:n2} which is a feasible flow for the PA formulation over node $i$. Similar to Example~\ref{exmp:n1}, we want to prevent such flows as Example~\ref{exmp:n2}.

\begin{example}\label{exmp:n2}
This is an example feasible solution for the PA formulation for the flow on incoming and outgoing arcs of node $i$ in Figure~\ref{fig:eq fig}: $x_{2i}^{\{k_3,k_4\}} = x_{i3}^{\{k_2\}} = 1$ and flow on other incoming and outgoing arcs of node $i$ equals to zero.
\end{example}

By having such modification on the partial aggregation network, it is enough to add flow conservation constraints for the added artificial nodes to tackle the issue raised in Section~\ref{sec:ineq-mip}. The modified network structure labels a commodity that arrives at a node by a disaggregated arc and thus excludes flow on an aggregated arc that does not include that specific commodity. Specifically, equations~\eqref{eq:flow-eq-1}-\eqref{eq:flow-eq-3} defines such flow conservation constraints. Here, we define the flow on artificial arcs by $z$ variables as constraints~\eqref{eq:z-eq-1} and \eqref{eq:z-eq-2}.
\begin{subequations}
\begin{align}
	\nonumber & \left(\sum_{j\in \mathcal{N}_i^+:D\cap\mathcal{D}_b^{ij}\neq\emptyset} x_{ij}^D - \sum_{j\in \mathcal{N}_i^-:D\cap\mathcal{D}_b^{ji}\neq\emptyset} x_{ji}^D\right) + \\& \left(\sum_{C\in \hat{\mathcal{T}}_b^i:D\cap C\neq\emptyset} z_{DC}^{ib}  -\sum_{C\in \check{\mathcal{T}}_b^i:D\cap C\neq\emptyset} 
    z_{CD}^{ib}\right) = \sum_{k \in D}(o_i^k-s_i^k)d^k  
	&&\forall b \in \mathcal{B},\ i\in\mathcal{C}_b,\ D \in \mathcal{M}_{b}^i:\mathcal{L}_b^i\neq\emptyset \label{eq:flow-eq-1}\\ 
	& \sum_{D\in \mathcal{M}_b^i:C\cap D\neq \emptyset}z_{CD}^{ib} - \sum_{j\in\mathcal{N}_i^-:C=\mathcal{K}_b^{ji}} x_{ji}^C = 0  &&\forall b \in \mathcal{B},\ i \in \mathcal{C}_b,\ C \in \check{\mathcal{T}}_{b}^i:\mathcal{L}_b^i\neq\emptyset \label{eq:flow-eq-2}\\
	& \sum_{j\in\mathcal{N}_i^+:C=\mathcal{K}_b^{ij}} x_{ij}^C    - \sum_{D\in \mathcal{M}_b^i:C\cap D\neq \emptyset}z_{DC}^{ib} = 0 &&\forall b \in \mathcal{B},\ i \in \mathcal{C}_b,\ C \in \hat{\mathcal{T}}_{b}^i:\mathcal{L}_b^i\neq\emptyset \label{eq:flow-eq-3}\\ 
	&z_{CD}^{ib} \geq 0 && \hspace{-4cm}\forall b\in \mathcal{B},\ i \in \mathcal{C}_b,\ C \in \check{\mathcal{T}}_{b}^i,\ D\in \mathcal{M}_b^i : C\cap D \neq \emptyset\ \&\ \mathcal{L}_b^i\neq \emptyset \label{eq:z-eq-1}\\
	&z_{DC}^{ib} \geq 0 && \hspace{-4cm}\forall b\in \mathcal{B},\ i \in \mathcal{C}_b,\ C \in \hat{\mathcal{T}}_{b}^i,\ D\in \mathcal{M}_b^i : C\cap D \neq \emptyset\ \&\ \mathcal{L}_b^i\neq \emptyset \label{eq:z-eq-2}
\end{align}
\end{subequations}

Equation~\eqref{eq:flow-eq-1} represents the flow conservation constraint for artificial intermediate nodes. The first part determines the net outgoing flow by the disaggregated arcs. The second part determines the net outgoing flow by the artificial arcs. The total outgoing flow must respect the corresponding demand at the node $i$ for the corresponding commodity set that an intermediate node represents. Equation~\eqref{eq:flow-eq-2} shows the flow conservation constraint for artificial in-flow aggregated nodes. The outgoing flow for such arcs is on artificial arcs and the incoming flow is on original aggregated arcs of the  partial aggregation network. Equation~\eqref{eq:flow-eq-3} shows the flow conservation constraint for artificial out-flow aggregated nodes. The outgoing flow for such arcs is on original aggregated arcs of the  partial aggregation network and the incoming flow is on artificial arcs. Decision variable $z$ is defined for the flow on an artificial arc where ever an artificial arc is added based on constraints \eqref{eq:z-eq-1} and \eqref{eq:z-eq-2}. Constraints~\eqref{eq:flow-eq-1}-\eqref{eq:z-eq-2} are added if  node $i$ of dispersion $b$ is required to be modified, i.e. $\mathcal{L}_b^i\neq \emptyset$. In the rest of the paper, we consider \eqref{eq:ag:obj}-\eqref{eq:ag:y_binary} and \eqref{eq:flow-eq-1}-\eqref{eq:z-eq-2} as the \emph{partially aggregated formulation with equality tightening  constraints}, and refer to it as \textbf{PAe}. Example~\ref{exmp:n2} is not feasible for Constraints~\eqref{eq:flow-eq-1}-\eqref{eq:z-eq-2} as those flow conservation constraints for the modified node $i$ do not link the flows $x_{2i}^{\{k_3,k_4\}}$ and~$x_{i3}^{\{k_2\}}$.

\subsection{K-shortest path aggregations}\label{sec:k-path}

The choices to cluster commodities can significantly affect the LP bound of a partially-aggregated formulation. In this section, we present a heuristic that employs a K-shortest path algorithm to effectively cluster the commodities. Based on Definition~\ref{def:cluster}, a dispersion aggregates a set of commodities, but it disaggregates each commodity from the group on a subset of arcs. Therefore, to construct a dispersion of a set of commodities, it is enough to determine a subset of arcs that each commodity of the considered set is disaggregated on. We call the subset of arcs on which the commodity $k$ is disaggregated from the corresponding group as \emph{critical arcs} and represent it by $\mathcal{A}^k$. Given $\mathcal{A}^k$ for all $k \in \mathcal{K}_b$, sets $\mathcal{K}_b^{ij}$ and $\mathcal{D}_b^{ij}$ of dispersion $b$ of the commodity set $\mathcal{K}_b$ are obtained as $D_b^{ij} = \{k\in \mathcal{K}_b\ \mid (i,j)\in \mathcal{A}^k\}$ and $\mathcal{K}_b^{ij} = \mathcal{K}_b \setminus D_b^{ij}$ for all $(i,j) \in \mathcal{A}$.

We determine the set of critical arcs of commodity $k$ by solving a K-shortest path problem from its origin $o^k$ to its destination $s^k$ considering the surrogate cost of $\Tilde{c}_{ij}=c_{ij} + \frac{f_{ij}}{u_{ij}}$ for all arcs $(i,j)\in\mathcal{A}$. Disaggregating commodity $k$ on arc $(i,j)$ translates to keeping the corresponding SI of the DA formulation in the partially-aggregated formulation. The surrogate cost $\Tilde{c}_{ij}$ incorporates the variable cost, fixed cost, and capacity of an arc simultaneously, and results in keeping impactful SIs that tighten the partially-aggregated formulation. Such surrogate cost has been also utilized in other solution algorithms for MCND. \cite{GENDRON201870} uses $\Tilde{c}_{ij}$ as the initial values for their slope scaling heuristic. \cite{Yaghini} uses $\Tilde{c}_{ij}$ to determine the set of variables that exist in the initial step of their column generation based heuristic. We have also examined two other cost structures as $c_{ij}+\frac{f_{ij}}{d^k}$ and $c_{ij}+\frac{f_{ij}}{\min_{k\in \mathcal{K}_b} d^k}$. However, $\Tilde{c}_{ij}$ appears to be significantly more promising, and therefore other cost structures are not considered in the rest of the paper. The details of the heuristic is explained in Algorithm~\ref{alg:k-path-agg}. This algorithm essentially considers the arcs that construct K shortest paths from origin to destination of the commodity $k$ as the set of its critical arcs $\mathcal{A}^k$, and forms the dispersions and the partial aggregation accordingly. We use the algorithm by \cite{kshortest} for finding K shortest paths between an origin and a destination. 

\begin{algorithm}[htb!]
\text{\textbf{Input: } $G=(\mathcal{N},\mathcal{A})$, set of commodities $\mathcal{K}$}, and the number of shortest paths considered $K$\\
\text{\textbf{Output: } A partial aggregation $\mathcal{B}$ which is a set of dispersions}\\
 \For{$b \in \mathcal{N}$} {
     $\mathcal{K}_b = \{k \in \mathcal{K} \mid o^k = b\}$\\
     \For{$k \in \mathcal{K}_b$}{
          $\mathcal{A}^k =$ Set of arcs included in K shortest paths from origin $o^k$ to destination $s^k$ considering the cost $\Tilde{c}_{ij}$ for each arc }
    \For {$(i,j) \in \mathcal{A}$}{
          $\mathcal{D}_b^{ij} = \{k\in \mathcal{K}_b\ \mid (i,j)\in \mathcal{A}^k\}$ \\ $\mathcal{K}_b^{ij} = \mathcal{K}_b \setminus \mathcal{D}_b^{ij}$}
    Add the dispersion $b$ to the set $\mathcal{B}$}
 \caption{K-shortest path aggregation}
 \label{alg:k-path-agg}
\end{algorithm}

\section{Polyhedral analysis}\label{sec:polyhed}
This section compares the polyhedra of the LP relaxations of the DA, PAe, PAi, PA, and FA formulations. We consider the LP relaxations of these formulation which are obtained by relaxing constraints~\eqref{eq:y_binary} and \eqref{eq:ag:y_binary} to $0\leq y_{ij} \leq 1$.
In the each of Theorems~\ref{thrm:da}-\ref{thrm:pa} we compare the strength of the LP relaxation of two formulations. To prove one of these two is stronger than the other, we show that there always exists a mapping that transforms an arbitrary feasible solution of the LP relaxation of the stronger formulation to a feasible solution for the LP relaxation of the weaker formulation. However, there exists fractional solutions to the weaker formulation that are not feasible for the stronger formulation. 

\begin{theorem}\label{thrm:da}
The DA formulation is stronger than the PAe formulaion.
\end{theorem}

\begin{proof}
Let $(\bar x_{ij}^k, \bar y_{ij})$ be an arbitrary solution of the LP relaxation of the DA formulation defined by~\eqref{eq:flow_conserv}-\eqref{eq:y_binary}. The following mapping transforms this solution to the solution ($x_{ij}^D, y_{ij}, z_{CD}^{ib}$) for the LP relaxation of PAe formulation defined by~\eqref{eq:ag:flow_conserv}-\eqref{eq:ag:y_binary} and \eqref{eq:flow-eq-1}-\eqref{eq:z-eq-2}.
\begin{subequations}
\begin{align}
    &x_{ij}^D \coloneqq  \sum_{k \in D} \bar x_{ij}^k &&\forall\ (i,j)\in \mathcal{A},\ b \in \mathcal{B},\ D \in \mathcal{G}_b^{ij} \label{eq:thm-da-x}\\
    &y_{ij} = \bar y_{ij} && \forall\ (i,j)\in \mathcal{A} \label{eq:thm-da-y}\\
    &z_{CD}^{ib} = \sum_{j\in \mathcal{N}^-_i}\sum_{k\in C\cap D} \bar x_{ji}^k && \forall b\in \mathcal{B},\ i \in \mathcal{C}_b,\ C \in \check{\mathcal{T}}_{b}^i,\ D\in \mathcal{M}_b^i : C\cap D \neq \emptyset\ \&\ \mathcal{L}_b^i\neq \emptyset  \label{eq:thm-da-z1}\\
    &z_{DC}^{ib} = \sum_{j\in \mathcal{N}^+_i}\sum_{k\in D\cap C} \bar x_{ij}^k  && \forall b\in \mathcal{B},\ i \in \mathcal{C}_b,\ C \in \hat{\mathcal{T}}_{b}^i,\ D\in \mathcal{M}_b^i : C\cap D \neq \emptyset\ \&\ \mathcal{L}_b^i\neq \emptyset \label{eq:thm-da-z2}
\end{align}
\end{subequations}
We now prove that the transformed solution ($x_{ij}^D, y_{ij}, z_{CD}^{ib}$) obtained by mapping~\eqref{eq:thm-da-x}-\eqref{eq:thm-da-z2} is feasible for LP relaxation of PAe. For a node related to the dispersion layer $b$, the flow conservation constraint~\eqref{eq:ag:flow_conserv} of PAe is just a summation of flow conversation constraints~\eqref{eq:flow_conserv} of DA over $k\in \mathcal{K}_b$.  The capacity constraint~\eqref{eq:ag:capacity} is feasible over the transformed solution since $\sum_{b\in \mathcal{B}}\sum_{D\in\mathcal{G}_b^{ij}} x_{ij}^D = \sum_{k\in \mathcal{K}} \bar x_{ij}^k \leq u_{ij} \bar y_{ij} = u_{ij} y_{ij} \ \forall (i,j) \in \mathcal{A}$. Strong inequality~\eqref{eq:ag:SI} is just a summation of SIs~\eqref{eq:stong_ineq} over $k\in D$, and therefore feasible for the transformed solution.

We then show that the transformed solution is feasible for constraints~\eqref{eq:flow-eq-1}-\eqref{eq:flow-eq-3}. Constraint~\eqref{eq:flow-eq-1} defines a flow conservation constraint for each artificial intermediate node (set $\mathcal{M}_b^i$) that is added to the network. Based on equations~\eqref{eq:thm-da-z1} and \eqref{eq:thm-da-z2}, Constraint~\eqref{eq:flow-eq-1} in fact conserves the flow of commodities of set $D$ at node $i$. Therefore, it is just a summation of Constraint~\eqref{eq:flow_conserv} over the commodities of the corresponding commodity set $D$. 

Since the set $\mathcal{M}_b^i$ partition the commodity set $\mathcal{K}_b$, in Constraint~\eqref{eq:flow-eq-2} we have:
$$\sum_{D\in \mathcal{M}_b^i:C\cap D\neq \emptyset}z_{CD}^{ib} = \sum_{j \in \mathcal{N}^-_i}\sum_{k \in C} \bar x^k_{ji}$$

On other hand, $\sum_{j\in\mathcal{N}_i^-:C=\mathcal{K}_b^{ji}} x_{ji}^C = \sum_{j \in \mathcal{N}^-_i}\sum_{k \in C} \bar x^k_{ji}$ based on equation~\eqref{eq:thm-da-x}. Therefore, Constraint~\eqref{eq:flow-eq-2} is satisfied by the transformed solution. In a similar way, Constraint~\eqref{eq:flow-eq-3} can be proved to be satisfied by the transformed solution.

Now we show an example where a solution of the PAe LP relaxation is not feasible for the DA LP relaxation. Consider an example that the commodity $k_1$ is aggregated with some other commodities into the set~$C$ on arc $(i,j)$. Moreover, $u_{ij} > \sum_{k \in C} d^k $. A feasible solution for PAe could be flowing only the demand of $k_1$ on arc $(i,j)$ that translates to $x_{ij}^C = d^{k_1}$. In this case, a feasible value for $y_{ij}$ is $\frac{d^{k_1}}{\sum_{k \in C} d^k}$ based on the corresponding SI and capacity constraint. However, this solution ($\bar x_{ij}^{k_1}= d^{k_1}, \bar y_{ij}=\frac{d^{k_1}}{\sum_{k \in C} d^k}$) is not feasible for the LP relaxation of DA as it violates the corresponding SI of the commodity $k_1$ on arc $(i,j)$.

Therefore, the DA formulation is stronger than the PAe formulation because of two reasons: 
first, there always exists a mapping that transforms an arbitrary solution of the LP relaxation of the DA formulation to a solution for the LP relaxation of the PAe formulation.  
Second, there exists feasible solutions for the LP relaxation of the PAe formulation that are not projectable to the LP relaxation of the DA formulation.
\end{proof}

\begin{theorem}\label{thrm:pae}
The PAe formulation is stronger than the PAi formulaion, where both formulations are based on a same partial aggregation $\mathcal{B}$.
\end{theorem}

\begin{proof}
Let $(\bar x_{ij}^D, \bar y_{ij}, \bar z_{CD}^{ib})$ be an arbitrary solution of the LP relaxation of the PAe formulation defined by~\eqref{eq:ag:flow_conserv}-\eqref{eq:ag:y_binary} and \eqref{eq:flow-eq-1}-\eqref{eq:z-eq-2}. A feasible solution ($x_{ij}^D, y_{ij}$) for the LP relaxation of PAi formulation defined by~\eqref{eq:ag:flow_conserv}-\eqref{eq:ag:y_binary}, \eqref{eq:path-ineq-1}, and \eqref{eq:path-ineq-2}  can be obtained by the mapping $x_{ij}^D \coloneqq \bar x_{ij}^D$ and $y_{ij} \coloneqq \bar y_{ij}$.

Constraints~\eqref{eq:ag:flow_conserv}-\eqref{eq:ag:y_binary} are common in both formulations. Therefore, we just prove that a solution by LP relaxation of PAe also satisfies constraints~\eqref{eq:path-ineq-1} and \eqref{eq:path-ineq-2}. Consider the commodity $k$ and its corresponding dispersion $b$ over node $i$ where $\{k\}\in\mathcal{M}_b^i$. Therefore, as in a partial aggregation the commodity $k$ on the arc $(i,j)$ is either aggregated with some other commodities or disaggregated, i.e. only one of two statements $(k \in \mathcal{K}_b^{ij} \Rightarrow \exists! D\in \mathcal{G}_b^{ij} : k\in  D)$ or $(k \in \mathcal{D}_b^{ij} \Rightarrow \exists! D\in \mathcal{G}_b^{ij} : D=\{k\})$ is correct for all $(i,j)\in \mathcal{A}$, we have:
\begin{subequations}
\begin{align}
    \sum_{j\in \mathcal{N}_i^+:D\cap\mathcal{D}_b^{ij}\neq\emptyset} x_{ij}^D = \sum_{j\in \mathcal{N}_i^+} \sum_{D\in \mathcal{G}_b^{ij}:D=\{k\}} x_{ij}^D \label{eq:thm:pae-x-diss-out}\\
    \sum_{j\in \mathcal{N}_i^-:D\cap\mathcal{D}_b^{ji}\neq\emptyset} x_{ji}^D = \sum_{j\in \mathcal{N}_i^-} \sum_{D\in \mathcal{G}_b^{ji}:D=\{k\}} x_{ji}^D \label{eq:thm:pae-x-diss-in}
\end{align}
Furthermore, based on equation~\eqref{eq:pae-out-agg-nodes} and considering $D=\{k\}$, we have $\{C\in\hat{\mathcal{T}}_b^i:D\cap C\neq \emptyset\} = \{\mathcal{K}_b^{ij} : \forall j \in\mathcal{N}_i^+ \mid k \in \mathcal{K}_b^{ij} \}$. Moreover, equation~\eqref{eq:flow-eq-3} implies that $z_{DC}^{ib} \leq \sum_{j\in\mathcal{N}_i^+:C=\mathcal{K}_b^{ij}} x_{ij}^C$ for all $C \in \hat{\mathcal{T}}_b^i$. Therefore, we have:
\begin{align}
    \sum_{C\in\hat{\mathcal{T}}_b^i:D\cap C\neq \emptyset} z_{DC}^{ib} \leq \sum_{j \in \mathcal{N}_i^+:k\in\mathcal{K}_b^{ij}} x_{ij}^{\mathcal{K}_b^{ij}} \label{eq:thm:pae-x-agg-out}
\end{align}
\end{subequations}
From Definition~\ref{def:dispersion-layer} we have $\mathcal{G}_b^{ij}=\mathcal{K}_b^{ij}\cup\{\{k\}\mid k \in \mathcal{D}_b^{ij}\}$. Moreover, $z_{CD}^{ib}\geq0$. Having these in mind and also statements~\eqref{eq:thm:pae-x-diss-out}-\eqref{eq:thm:pae-x-agg-out}, flow conservation constraint~\eqref{eq:flow-eq-1} of PAe implies the inequality~\eqref{eq:path-ineq-1} of PAi. Therefore, any arbitrary solution of PAe satisfies Constraint~\eqref{eq:path-ineq-1} of the PAi formulation. Similar approach can prove that the feasibility of solutions of PAe for Constraint~\eqref{eq:path-ineq-2} of PAi. Clearly, when $\{k\}\notin\mathcal{M}_b^i$, constraints~\eqref{eq:path-ineq-1} and \eqref{eq:path-ineq-2} are redundant regarding flow conservation constraint~\eqref{eq:ag:flow_conserv}. 

Now we show an example where a solution of the PAi LP relaxation that is not feasible for the PAe LP relaxation. Consider an example that dispersion $b$ of the commodity set $\mathcal{K}_b=\{k_1,k_2,k_3,k_4\}$ over node $i$, where $\mathcal{N}_i^-=\{1,2\}$ and $\mathcal{N}_i^+=\{3,4\}$. All commodities of $\mathcal{K}_b$ are disaggregated on arcs $(1,i)$ and $(i,3)$, i.e. $\mathcal{K}_b^{1i}=\mathcal{K}_b^{i3}=\emptyset$ and $\mathcal{D}_b^{1i}=\mathcal{D}_b^{i3}=\mathcal{K}_b$. However, they are all aggregated on two other arcs, i.e. $\mathcal{K}_b^{2i}=\mathcal{K}_b^{i4}=\mathcal{K}_b$ and $\mathcal{D}_b^{2i}=\mathcal{D}_b^{i4}=\emptyset$. A feasible solution by PAi for this example is ($x_{1i}^{\{k_1\}}=x_{1i}^{\{k_2\}}=1$, $x_{1i}^{\{k_3\}}=x_{1i}^{\{k_4\}}=0$, $x_{i3}^{\{k_1\}}=x_{i3}^{\{k_2\}}=0$, $x_{i3}^{\{k_3\}}=x_{i3}^{\{k_4\}}=1$, $x_{2i}^{\mathcal{K}_b}=x_{i4}^{\mathcal{K}_b}=1$). However, this solution is infeasible for the additional flow conservation constraints~\eqref{eq:flow-eq-1}-\eqref{eq:flow-eq-3} of the PAe formulation as shown in the modified node $i$ in  Figure~\ref{fig:proof-example}. In this figure, the flow on each is arc is labeled on it inside the parenthesis. Considering a flow conservation constraint for each internal node, the solution shown on the arcs is not feasible for PAe LP relaxation. 
\begin{figure}[htb!]
    	\centering
        \begin{tikzpicture}
        \draw (-4,-6.75) rectangle (4, 0.75);
        \node[] at (0,1) () {};
        \node[circle,draw,thick,dashed,blue,minimum size = 1cm] at (0,0) (m-1) {\tiny $k_1$};
        \node[circle,draw,thick,dashed,red,minimum size = 1cm] at (0,-1.5) (m-2) {\tiny $k_2$};
        \node[circle,draw,thick,dashed,Green4,minimum size = 1cm] at (0,-3) (m-3) {\tiny $k_3$};
        \node[circle,draw,thick,dashed,DarkOrchid4,minimum size = 1cm] at (0,-4.5) (m-4) {\tiny $k_4$};
        \node[circle,draw,thick,dashed,minimum size = 1cm] at (-3,-6) (i-1) {\tiny $\mathcal{K}_b$};
        \node[circle,draw,thick,dashed,minimum size = 1cm] at (3,-6) (o-1) {\tiny $\mathcal{K}_b$};
        \path[->, thick,dashed]   (i-1)  edge (m-1);
        \path[->, thick,dashed]   (i-1)  edge (m-2);
        \path[->, thick,dashed]   (i-1)  edge (m-3);
        \path[->, thick,dashed]   (i-1)  edge (m-4);
        \path[->, thick,dashed]   (m-1)  edge (o-1);
        \path[->, thick,dashed]   (m-2)  edge (o-1);
        \path[->, thick,dashed]   (m-3)  edge (o-1);
        \path[->, thick,dashed]   (m-4)  edge (o-1);
        \node[circle,draw,thick,minimum size = 1cm] at (-6,-2) (n1) {$1$};
        \node[circle,draw,thick,minimum size = 1cm] at (-6,-4) (n2) {$2$};
        \node[circle,draw,thick,minimum size = 1cm] at (6,-2) (n3) {$3$};
        \node[circle,draw,thick,minimum size = 1cm] at (6,-4) (n4) {$4$};
        \path[->, thick,blue]   (n1)  edge  node[above,sloped]{\footnotesize$\{k_1\}, (1)$}(m-1);
        \path[->, thick,red]   (n1)  edge  node[above,sloped]{\footnotesize$\{k_2\}, (1)$}(m-2);
        \path[->, thick,Green4]   (n1)  edge  node[above,sloped]{\footnotesize$\{k_3\}, (0)$}(m-3);
        \path[->, thick,DarkOrchid4]   (n1)  edge  node[above,sloped]{\footnotesize$\{k_4\}, (0)$}(m-4);
        \path[->, thick,blue]   (m-1)  edge  node[above,sloped]{\footnotesize$\{k_1\}, (0)$}(n3);
        \path[->, thick,red]   (m-2)  edge  node[above,sloped]{\footnotesize$\{k_2\}, (0)$}(n3);
        \path[->, thick,Green4]   (m-3)  edge  node[above,sloped]{\footnotesize$\{k_3\}, (1)$}(n3);
        \path[->, thick,DarkOrchid4]   (m-4)  edge  node[above,sloped]{\footnotesize$\{k_4\}, (1)$}(n3);
        \path[->, thick]   (n2)  edge  node[above,sloped]{\footnotesize$\mathcal{K}_b, (1)$}(i-1);
        \path[->, thick]   (o-1)  edge  node[above,sloped]{\footnotesize$\mathcal{K}_b, (1)$}(n4);

		\end{tikzpicture}
        \caption{An example solution of PAi LP relaxation that is not feasible for PAe LP relaxation. This figure shows the modified node $i$ for PAe, where $\mathcal{N}_i^-=\{1,2\}$ and $\mathcal{N}_i^+=\{3,4\}$. The dispersion $b$ (including commodity set $\mathcal{K}_b=\{k_1,k_2,k_3,k_4\}$) of a partial aggregation states that $\mathcal{K}_b^{1i}=\mathcal{K}_b^{i3}=\emptyset$, $\mathcal{D}_b^{1i}=\mathcal{D}_b^{i3}=\mathcal{K}_b$, $\mathcal{K}_b^{2i}=\mathcal{K}_b^{i4}=\mathcal{K}_b$, and $\mathcal{D}_b^{2i}=\mathcal{D}_b^{i4}=\emptyset$. A feasible solution of the PAi LP relaxation for the flow of arcs is shown as the numbers in parenthesis above each arc. However, this solution is infeasible for PAe because of the flow conservation constrains for the internal nodes.} 
        \label{fig:proof-example}
    \end{figure}
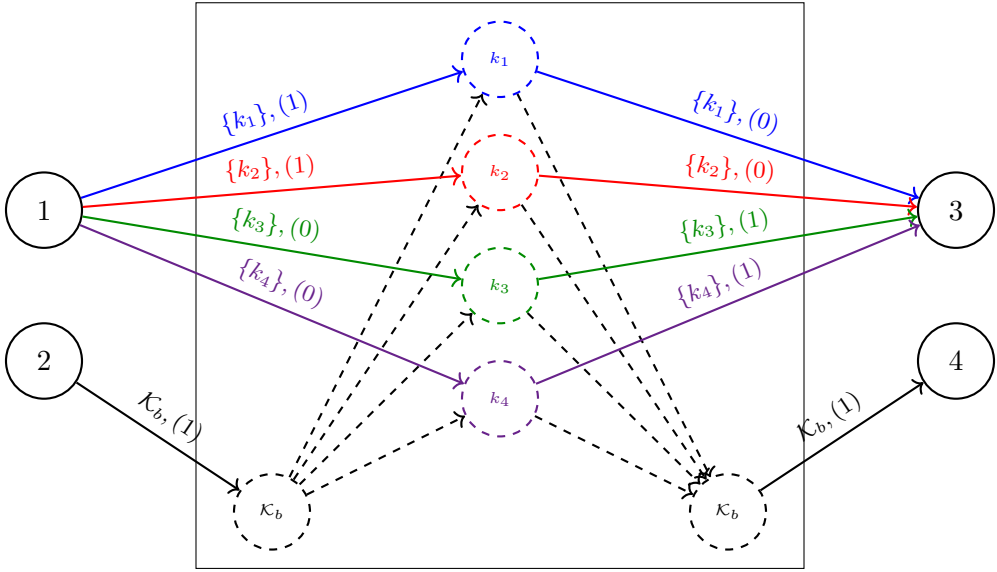
Therefore, the PAe formulation is stronger than the PAi formulation.
\end{proof}

\begin{theorem}\label{thrm:pai}
The PAi formulation is stronger than the PA formulation, where both formulations are based on a same partial aggregation $\mathcal{B}$.
\end{theorem}

\begin{proof}
PA formulation could be obtained from PAi by dropping constraints\eqref{eq:path-ineq-1} and \eqref{eq:path-ineq-2}. Therefore, any feasible solution of PAi is also feasible for PA formulation.

Now we show an example where a solution of the PA LP relaxation that is not feasible for the PAi LP relaxation. Consider node $i$ in Figure~\ref{fig:ineq fig b} and its corresponding partial aggregation. A feasible solution for PA LP relaxation over this partial aggregation over node $i$ is $x_{1i}^{\mathcal{K}_b\setminus \{k_1\}} = d^{k_1}$ and $x_{i3}^{\{k_1\}} =  d^{k_1}$. However, this solution violates the corresponding Constraint~\ref{eq:path-ineq-1} in the PAi LP relaxation.

Therefore, the PAi formulation is stronger than the PA formulation.
\end{proof}

\begin{theorem}\label{thrm:pa}
The PA formulation is stronger than the FA formulation.
\end{theorem}

\begin{proof}
Let $(\bar x_{ij}^D, \bar y_{ij})$ be an arbitrary solution of the LP relaxation of the PA formulation for an arbitrary partial aggregation $\mathcal{B}^p$. The following mapping transforms this solution to the solution ($x_{ij}^D, y_{ij}$) for the LP relaxation of FA formulation based on a full aggregation $\mathcal{B}^f=\Tilde{\mathcal{N}}$ as explained in Section~\ref{sec:commodity}. 
\begin{subequations}
\begin{align}
    &x_{ij}^{\mathcal{K}_n} \coloneqq  \sum_{b \in \mathcal{B}^p:o_b=n} \sum_{D\in \mathcal{G}_b^{ij}} \bar x_{ij}^D &&\forall\ (i,j)\in \mathcal{A},\ n \in \mathcal{B}^f \label{eq:thm-f-x}\\
    &y_{ij} = \bar y_{ij} && \forall\ (i,j)\in \mathcal{A} \label{eq:thm-f-y}
    \end{align}
\end{subequations}
The solution obtained by mappings \eqref{eq:thm-f-x} and \eqref{eq:thm-f-y} is feasible for the FA formulation. The flow conservation constraint of FA is a summation of flow conservation constraints of dispersion with the same origin. The capacity constraint in both cases is the same since the flow of all commodities are included in one inequality for each arc. Similar to the flow conservation constraints, SIs of the FA formulation are summation of the SIs of relevant dispersions with the same origin. Hence, an arbitrary solution of PA formulation is always feasible for the FA formulation as the FA formulation is obtained by aggregating some of the PA constraints.

Now we show an example where a solution of the DA LP relaxation that is not feasible for the PA LP relaxation. Consider an example of a partial aggregation that two commodities $k_1$ and $k_2$ are not aggregated on the arc $(i,j)$. Moreover, $d^{k_1}+d^{k_2}>u_{ij}$, but $d^{k_1}<u_{ij}$ and $d^{k_2}<u_{ij}$. As in the FA formulation these commodities are aggregated together, a feasible solution is $x_{ij}^{\{k_1,k_2\}} = d^{k_1}$ and $y_{ij}=\frac{d^{k_1}}{u_{ij}}$. However, this solution is not feasible for the PA formulation based on its SIs for the arc $(i,j)$.

Therefore, the PA formulation is stronger than the FA formulation.
\end{proof}

Theorems~\ref{thrm:da}-\ref{thrm:pa} imply a hierarchy of formulations in terms of LP relaxation strength. In this hierarchy, the DA formulation is the strongest formulation and the FA formulation is the weakest. Partially-aggregated formulations lie in between these two extreme cases. Moreover, as each commodity is included in one unique dispersion, we have $\sum_{b \in \mathcal{B}}\sum_{D\in\mathcal{G}_b^{ij}}x_{ij}^D=\sum_{k\in\mathcal{K}}x_{ij}^k$ which shows that the objective function of the partially-aggregated formulations is equivalent to the objective function of the DA formulation. Therefore, we can conclude Corollary~\ref{crll:valid-lb}.
 
\begin{corollary}\label{crll:valid-lb}
PA, PAi, and PAe formulations are valid bounds for MCND.
\end{corollary}

\section{Computational results}\label{sec:comp}

In this section we perform two sets of computational experiments to evaluate the proposed formulations. The first set of experiments investigates the LP relaxations of these formulations in terms of tightness and computing time. The second set of experiments compares the mixed-integer programming algorithms over the formulations. All computational experiments are run on a cluster with 4 Xeon-Gold-6150 cores and 8 GB RAM using CPLEX 12.10.0 via a Python API as LP and MIP solver. Default settings are always selected unless otherwise specified.

Computational experiments are conducted on 196 publicly available instances for MCND, called \emph{Canad} instances. These instances were generated by \cite{CRAINIC200173} and used as benchmark instances by subsequent papers to evaluate different solution algorithms and relaxations. These instances cover a diverse range of number of commodities,  capacity tightness, network density, and significance of fixed cost over the flow cost. The details of the instances can be found in \cite{CRAINIC200173}. These instances are divided into three classes. Class \emph{C} includes 31 large instances with dense networks and many commodities. Class \emph{C+} includes 12 instances with sparser networks and few commodities. Class \emph{R} includes 153 small and medium size instances with dense networks but a smaller ratio of the number of commodities per  number of nodes in comparison with the C class.

Table~\ref{tab:ins-summary} gives a summary of instances considered in this paper. In three instances of the C+ class, each node is the origin of at most one commodity. Therefore, we ignore these three instances as all formulations are the same for such cases. As the formulations mainly affect the computing time of the LP relaxations, we distinguish a set of instances as \emph{long} instances. A long instance is an instance that the computing time of its DA LP relaxation is longer than 10 seconds by the solver. 

\begin{table}[htb!]
\centering
\caption{Summary of the considered instances}
\label{tab:ins-summary}
\begin{tabular}{cccc}
\hline
\textbf{Class}       & \textbf{C}          & \textbf{C+} & \textbf{R}               \\ \hline
\#instances          & 31                  & 9           & 153                      \\
\#long instances     & 11                  & 0           & 7                        \\
\#nodes              & \{20,30\}           & \{25,100\}  & \{10,20\}                \\
\#arcs               & \{230,300,520,700\} & \{100,400\} & \{35,60,85,120,220,320\} \\
\#commodities        & \{40,100,200,400\}  & \{10,30\}   & \{10,25,40,50,100,200\}  \\ \hline
Arc density          & 70\%                & 11\%        & 63\%                     \\
Commodity density    & 31\%                & 2\%         & 31\%                     \\
Commodity/node ratio & 7.3                 & 0.6         & 4.3                      \\
FA reduction    & 78\%                & 24\%        & 67\%                     \\ \hline
\end{tabular}
\end{table}

The first part of Table~\ref{tab:ins-summary} shows the number of total and long instances of each class. It also describes the possible values for the number of nodes, arcs, and commodities of instances in each class. The second part reports the average value of some characteristics of instances in each group. \emph{Arc density} expresses how many of the possible arcs exist in an instance as $\frac{|\mathcal{A}|}{|\mathcal{N}|(|\mathcal{N}|-1)}$. \emph{Commodity density} presents the same notion for the commodities as $\frac{|\mathcal{K}|}{|\mathcal{N}|(|\mathcal{N}|-1)}$. The \emph{Commodity/node ratio} is the average number of commodities that originate from a node. \vphantom{$\frac{K}{N}$}
\emph{FA reduction} states the average reduction of the number of commodities in percentage across all instances by transforming the DA formulation ($|\mathcal{K}|$ commodities) to the FA formulation ($|\Tilde{\mathcal{N}}|$ commodities). FA reduction gives an indication of possible size reduction for an instance by a commodity aggregation scheme. As the numbers in Table~\ref{tab:ins-summary} suggest, there is more opportunity on average to reduce the size of the instances of C and C+ classes by commodity aggregations. 

As we use Algorithm~\ref{alg:k-path-agg} in the rest of the paper to obtain partially-aggregated formulations, we refer to PAi and PAe as PAi-K and PAe-K to develop and compare different formulations with respect to the number of shortest paths ($K$) used to construct dispersion of commodities.

\subsection{Experimental evaluation of the LP relaxations}\label{sec:lp-result}

In this section we evaluate the LP relaxation of the proposed formulations. LP relaxations are obtained by relaxing Constraint~\eqref{eq:y_binary} to $0\leq y_{ij} \leq 1$. First, we investigate the LP relaxations in terms of solution time, size, and LP bound by solving the LP model by CPLEX LP solver. Afterwards, the performance of the CPLEX cutting plane algorithm at the root node over the formulations is studied. In both cases, results indicate that formulations based on partial aggregation are on a Pareto frontier for the trade-off between the LP bound strength and the LP solution time. In our tests, the DA and FA formulation, which can be seen as extreme versions of both PAi-K and PAe-K, were compared against various forms of partial aggregation. These partial aggregations were created using the heuristic Algorithm~\ref{alg:k-path-agg} using the $K$ shortest paths solutions for $K\in\{1,\ldots,5,10\}$, giving rise to 12 partially aggregated formulations PAe-1,$\ldots$,PAe-10 and PAi-1,$\ldots$,PAi-10.

\begin{figure}[htb!]
    \centering
    \includegraphics[width=0.7\textwidth]{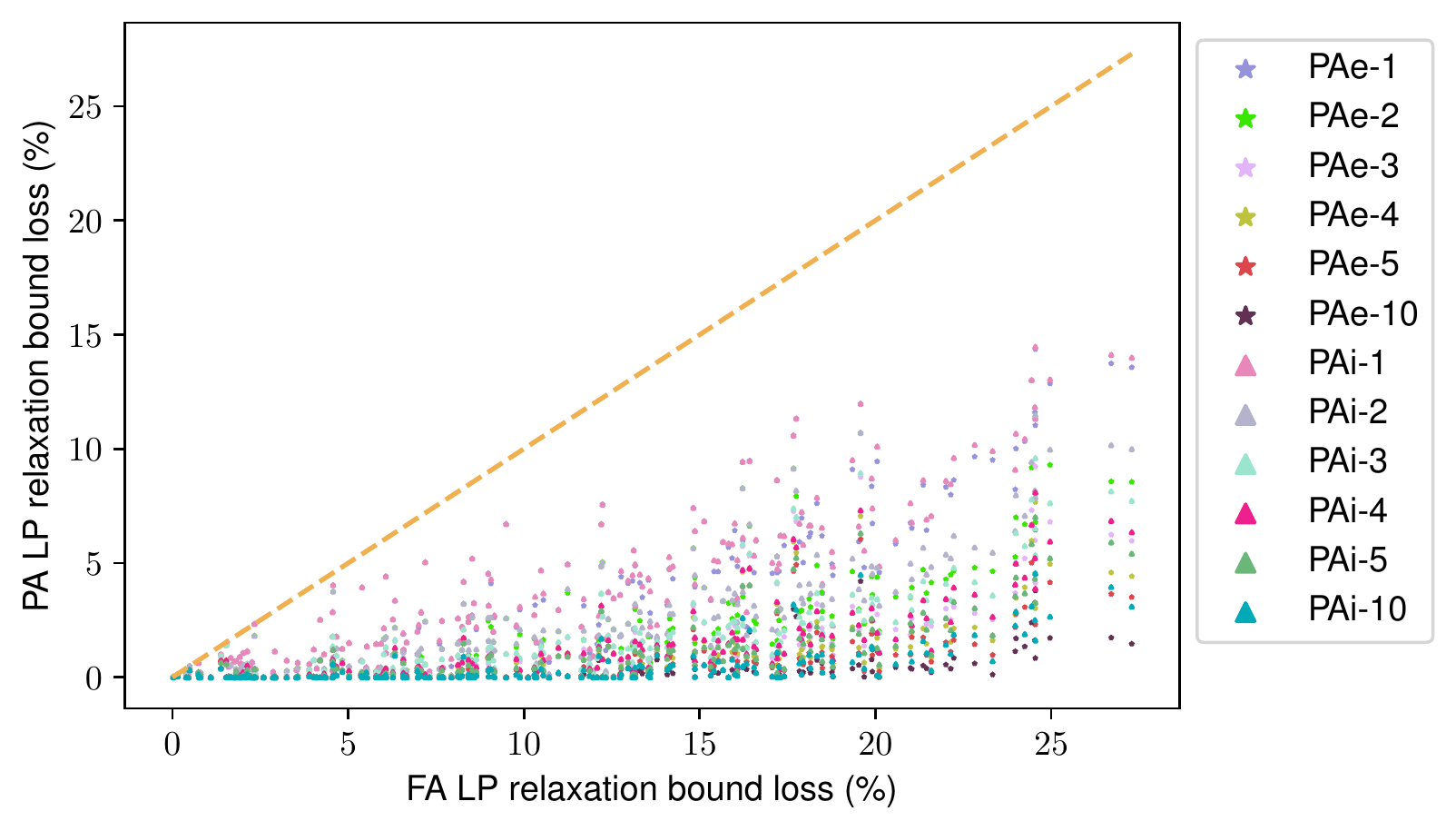}
    \caption{LP bound loss of partially-aggregated formulations versus LP bound loss of the FA formulation, considering the DA formulation as the base}
    \label{fig:lp-VSaggregation}
\end{figure}

\begin{figure}[htb!]
    \centering
    \includegraphics[width=0.7\textwidth]{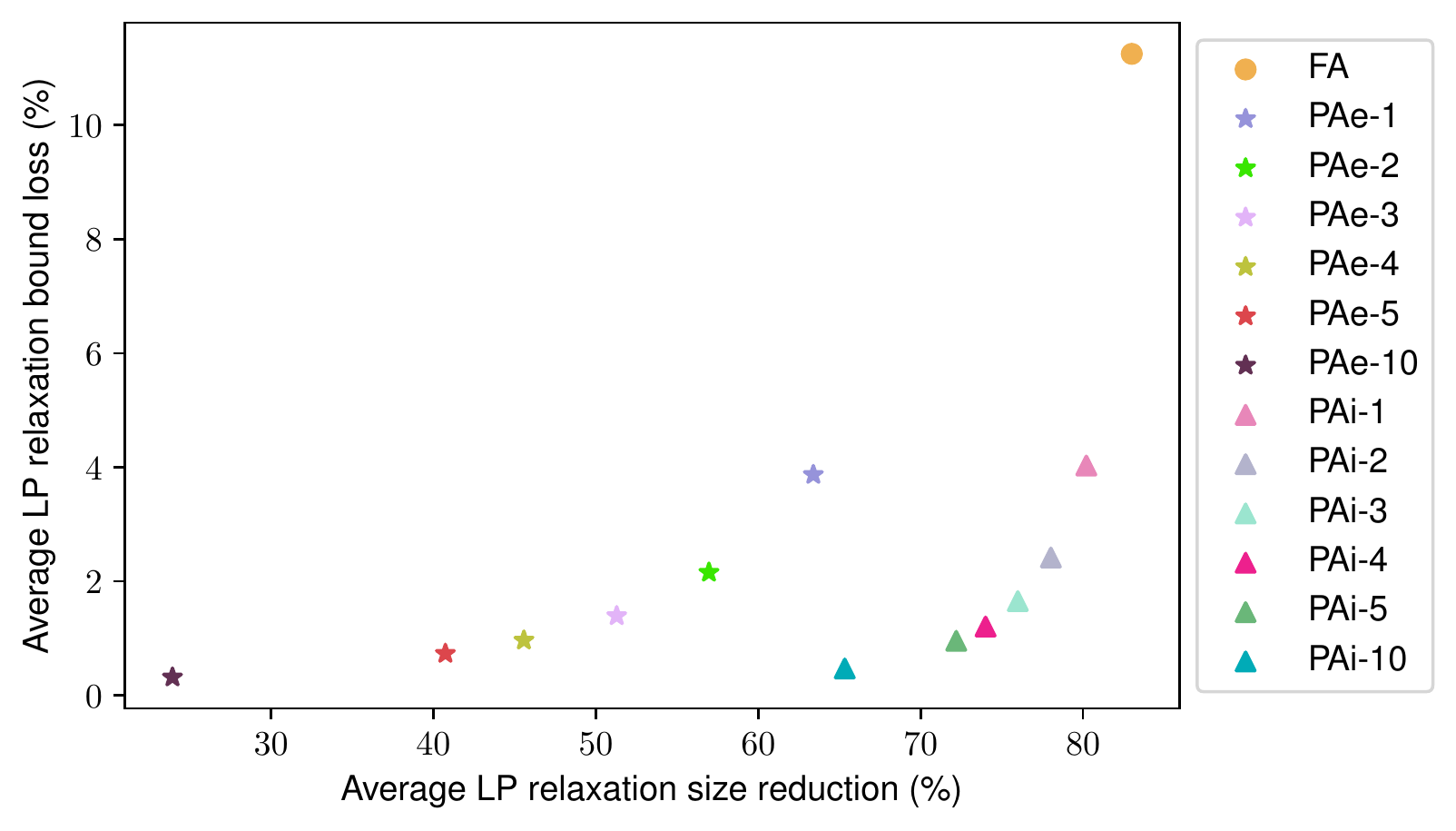}
    \caption{The trade-off between the LP bound loss and size for aggregated formulations considering the DA formulation as the base. Size of a formulation is considered as the multiplication of the number of variables and the number of constraints it includes.}
    \label{fig:lp-size}
\end{figure}

\begin{figure}[htb!]
    \centering
    \includegraphics[width=0.7\textwidth]{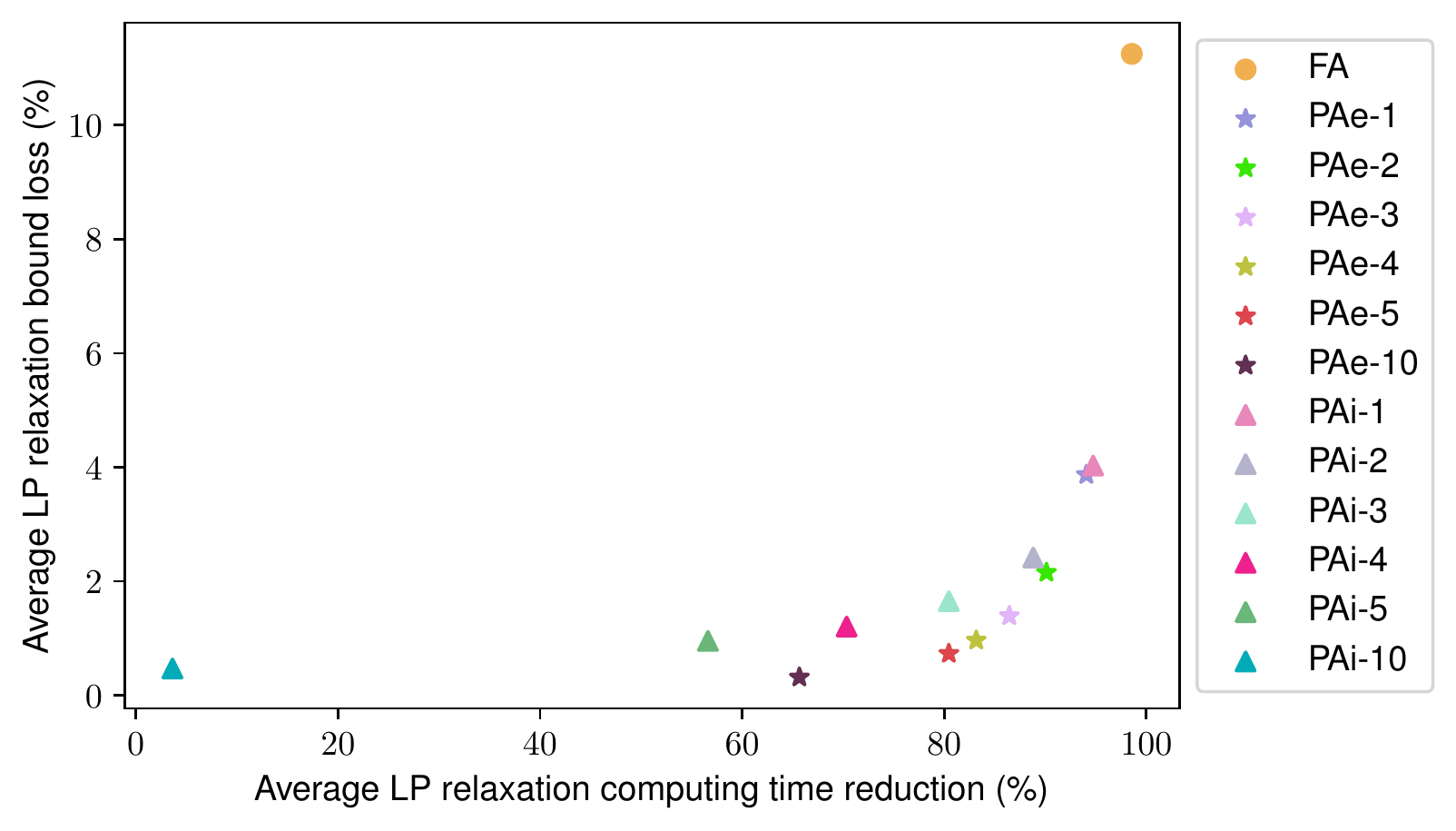}
    \caption{The trade-off between the LP bound loss and computing time for aggregated formulations considering the DA formulation as the base}
    \label{fig:lp-time}
\end{figure}

Figures~\ref{fig:lp-VSaggregation}-\ref{fig:lp-time} presents the results for solving the LP relaxations using CPLEX with default settings. In these figures, the LP bound of fully- and partially-aggregated formulations are compared with the LP bound of the DA formulation. 

\emph{Bound loss} measures the percentage of the LP optimal value decrease in a formulation with respect to the LP optimal value of the DA formulation. Figure~\ref{fig:lp-VSaggregation} compares the bound loss of the partially-aggregated formulations against the bound loss of the fully-aggregated formulation over all instances. The y-axis value of a point shows the bound loss of the corresponding formulation over an instance, whereas the x-axis value of this point represents the bound loss of the FA formulation over the corresponding instance. The results demonstrate the effectiveness of partially-aggregation formulations as their bound loss is significantly less than the bound loss of the FA formulation, even when considering just one shortest path to construct the partial aggregations.

Figures~\ref{fig:lp-size} and~\ref{fig:lp-time} indicate the trade-off between the bound loss and the LP size/computing time. \emph{Size reduction} is quantified over both dimensions, rows and columns, as the percentage change in the multiplication of the number of variables and the number of constraints of a formulation in comparison with the DA formulation. Similarly, \emph{time reduction} is with respect to the  computing time of the DA LP relaxation. These figures indicate that the size and solution time reduction of the partially-aggregated formulations are a substantial fraction of those achieved by the FA formulation. The time reduction is particularly significant over long instances. The average computing time of DA LP relaxation over long instances is 61.71 seconds. For instance, using the PAe-5 formulation reduces the average computing time by 85.6\% to 8.92 seconds in the expense of 1.2\% bound loss on average. On the other hand, the FA formulation reduces the computing time by 99.5\% but with a significant bound loss (17.3\% on average).

\begin{figure}[htb!]
    \centering
    \includegraphics[width=0.6\textwidth]{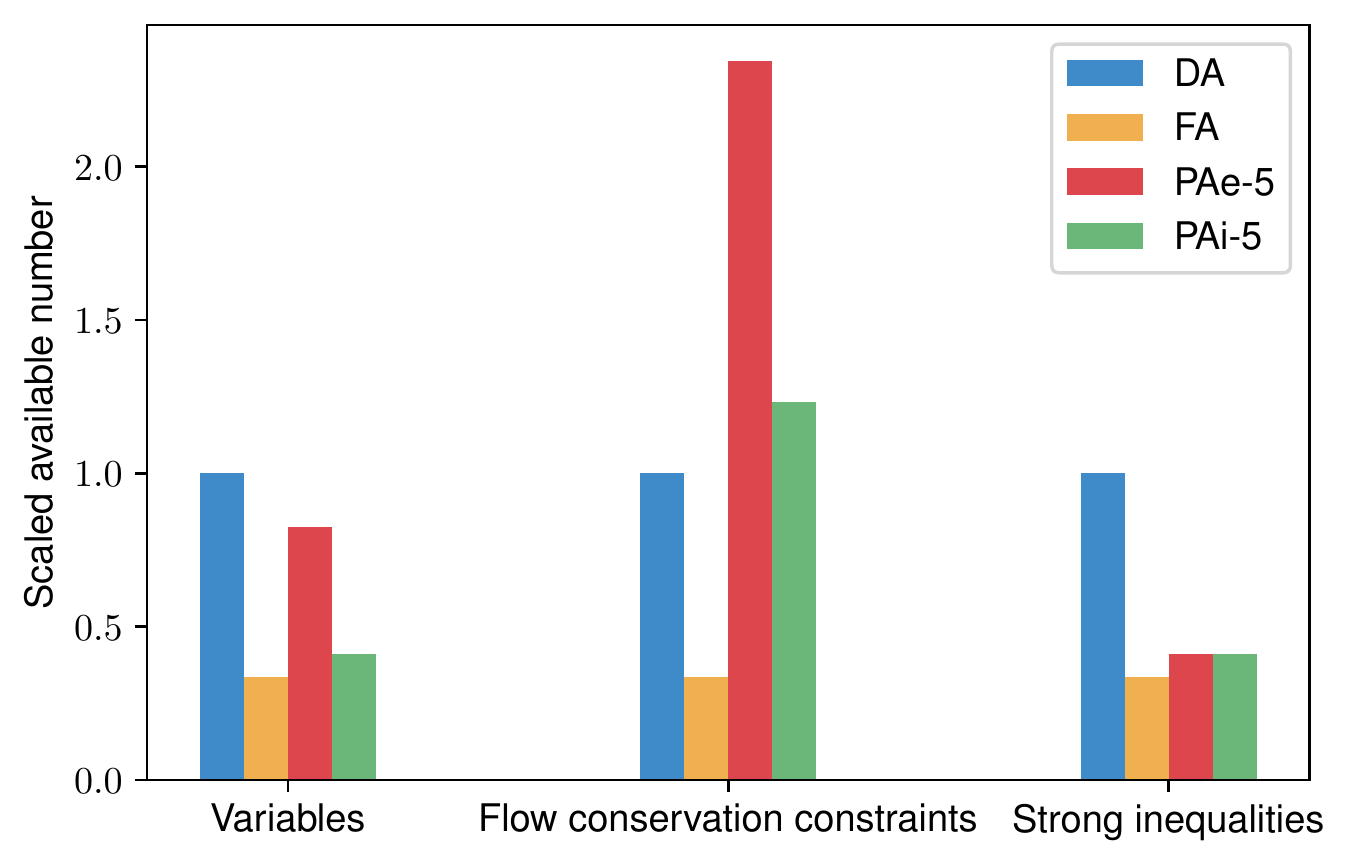}
    \caption{Structure of different formulations in terms of scaled number of variables and constraints, considering the dimensions of the DA formulation as the base}
    \label{fig:lp-structure}
\end{figure}

Figure~\ref{fig:lp-structure} depicts the average scaled number of variables, flow conservation constraints, and strong inequalities for all formulations considering the dimensions of the DA formulation as the base. Only the partial formulations constructed using 5 shortest paths are shown in this figure. However, the trend can be inferred for other values. Moreover, Constraints~\eqref{eq:path-ineq-1} and~\eqref{eq:path-ineq-2} are considered as flow conservation constraints for PAi-5. Partial formulations include more flow conservation constraints, particularly PAe-5. However, the number of strong inequalities and the number of variables are reduced leading to faster LP solve times. Mostly, strong inequalities make a difference in a formulation in terms of the LB bound and computing time. In fact, without any SI, the LP bound of all formulations are the same, and their computing time are in the same order. The number of SIs in the partial formulations (with the ratio of 0.41) is slightly more than the number of SIs in the FA formulation (with the ratio of 0.34). This implies an average of 59\% reduction in the number of SIs by partial formulations and 66\% reduction by the FA formulation. However, the slightly more SIs leads to significant error bound reduction as shown earlier. 

\begin{figure}[htb!]
    \centering
    \includegraphics[width=0.7\textwidth]{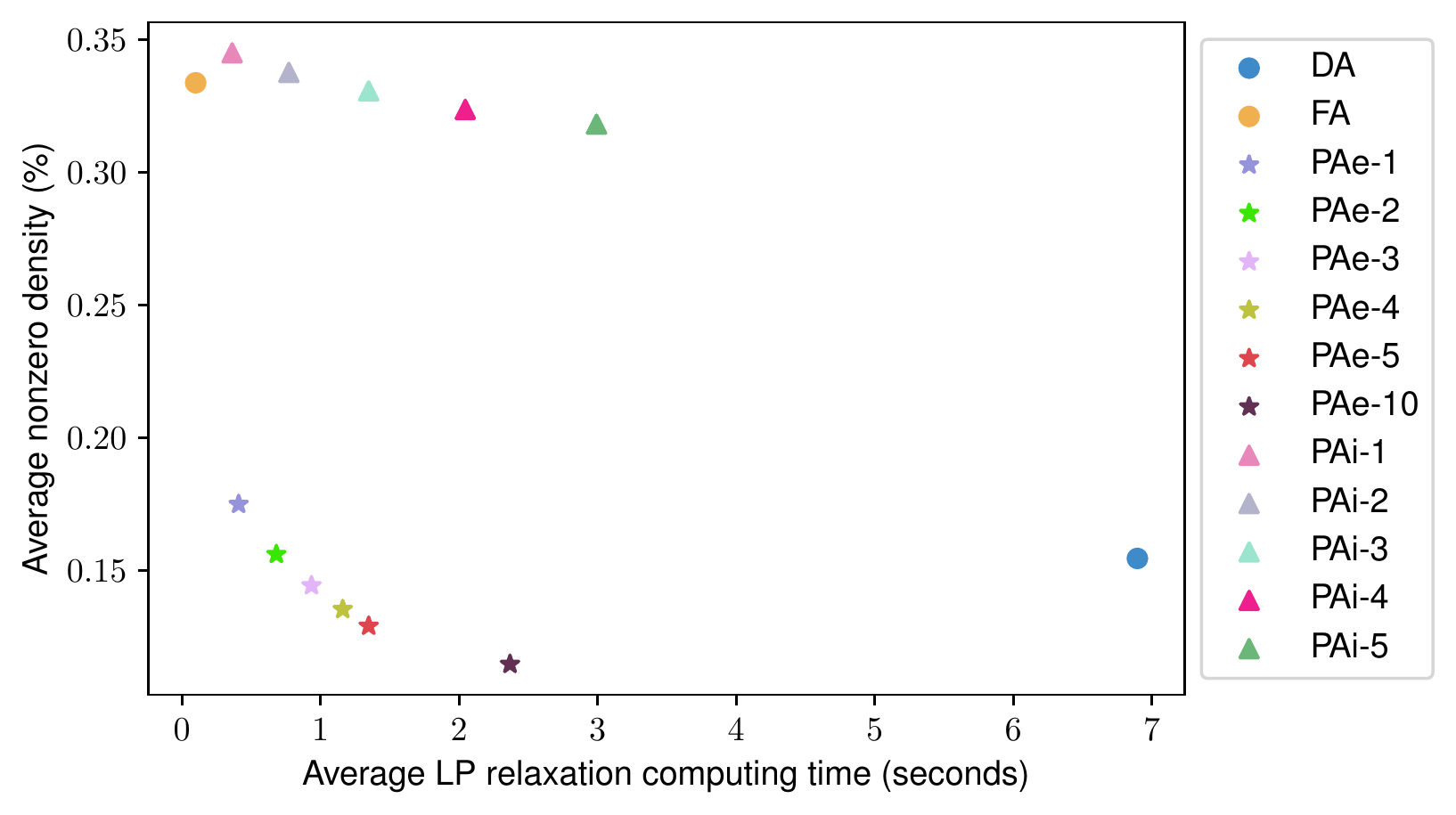}
    \caption{LP nonzeros density versus LP solution time}
    \label{fig:lp-nonzero}
\end{figure}

Because of the artificial nodes and arcs, the size of the PAe-K formulation is generally larger than the PAi-K formulation. However, the larger size does not translate to longer computing time as seen in Figure~\ref{fig:lp-time}. The shorter computing time of PAe-K despite a larger size is due to different structures of PAe-K and PAi-K. The PAe-K formulation, in contrast to PAi-K, keeps the multicommodity network flow structure  as it adds only equality flow conservation constraint. Another characteristic that affects the LP solution times is the density of nonzeros in the formulations. Figure~\ref{fig:lp-nonzero} illustrates the average nonzero density of each formulation with respect to the average solution time over all instances. Nonzero density is obtained by dividing the number of nonzeros by the size of the LP as the number of constraints multiplied by the number of constraints. In general, the nonzero density of PAe-K is less than the nonzero density of PAi-K which benefits the LP solution time. 

\begin{figure}[htb!]
    \centering
    \includegraphics[width=0.9\textwidth]{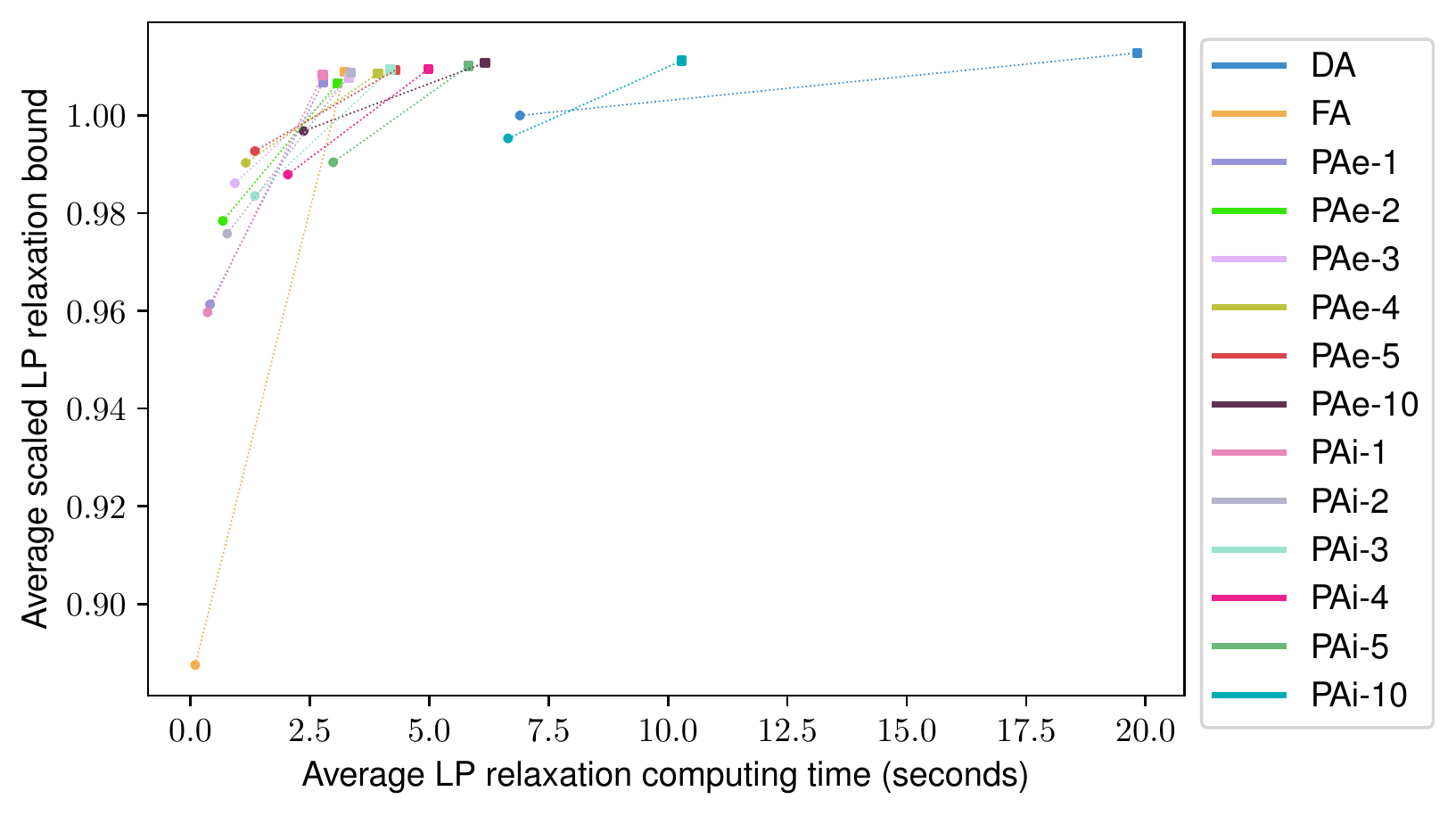}
    \caption{Average scaled LP bounds of the formulations with respect to their average computing time. Circle nodes correspond to the LP relaxation, and square nodes correspond to the LP bound by the cutting plane algorithm. The DA LP relaxation is considered as the base for scaling.}
    \label{fig:lp-pareto}
\end{figure}

The strength of the lower bound is significantly impacted by the addition of cuts during the MIP solution process. To assess the effect of this, Table~\ref{tab:cutting-plane} presents the summary results of the CPLEX cutting plane algorithm over all formulations. The results are obtained by restricting the CPLEX to 0 node of B\&B tree. Moreover, the best-known solutions are given as the initial incumbent to reduce the effect of heuristics on the computing times. Bound losses are with respect to the  DA LP relaxation. Therefore, negative bound losses in the  cutting plane algorithm column implies that the LP bound is improved in comparison with the  DA LP relaxation. Based on Table~\ref{tab:cutting-plane}, the significant bound loss of FA LP relaxation is compensated by the cutting plane algorithm. The DA formulation gives the best bound by using the cutting plane algorithm, however, by a remarkably higher computational effort. It is interesting that the partially- and fully-aggregated formulations, by using the cutting plane algorithm, are able to provide better bounds in comparison with the  DA LP relaxation in almost half computing time. Figure~\ref{fig:lp-pareto} displays the data from Table~\ref{tab:cutting-plane} in a graphical form by showing the average scaled LP bound for each formulation, considering the LP relaxation of DA formulation as the base. As Figure~\ref{fig:lp-pareto} shows that in both cases, the  LP bound and the LP bound by the cutting plane algorithm, there are partial formulations that are on the Pareto frontier for the trade-off between the LP bound and the computing time. The cutting plane algorithm is significantly faster over aggregated formulations in comparison with the DA formulation.

\begin{table}[htb!]
\caption{Effect of formulations on the performance of the CPLEX cutting plane algorithm}
\label{tab:cutting-plane}
\centering
\footnotesize
\begin{tabular}{c|cc|ccc} 
\hline
\multicolumn{1}{c|}{\multirow{2}{*}{Formulation}} & \multicolumn{2}{c|}{LP relaxation} & \multicolumn{3}{c}{Cutting plane   algorithm} \\ \cline{2-6} 
\multicolumn{1}{c|}{}    & Time  & Bound loss    & Time & Bound loss  & \#Cuts    \\ \hline
DA     & 6.9  & 0.00\%  & 19.8 & -1.28\%  & 56  \\
FA     & 0.1  & 11.25\% & 3.2  & -0.89\%  & 289 \\ \hline
PAe-1  & 0.4  & 3.87\%  & 2.8  & -0.68\%  & 355 \\
PAe-2  & 0.7  & 2.16\%  & 3.1  & -0.66\%  & 279 \\
PAe-3  & 0.9  & 1.39\%  & 3.3  & -0.76\%  & 230 \\
PAe-4  & 1.2  & 0.97\%  & 3.9  & -0.86\%  & 202 \\
PAe-5  & 1.3  & 0.73\%  & 4.3  & -0.93\%  & 177 \\
PAe-10 & 2.4  & 0.32\%  & 6.2  & -1.08\%  & 120 \\ \hline
PAi-1  & 0.4  & 4.03\%  & 2.8  & -0.83\%  & 242 \\
PAi-2  & 0.8  & 2.42\%  & 3.4  & -0.88\%  & 213 \\
PAi-3  & 1.3  & 1.65\%  & 4.2  & -0.94\%  & 194 \\
PAi-4  & 2.0  & 1.21\%  & 5.0  & -0.95\%  & 178 \\
PAi-5  & 3.0  & 0.96\%  & 5.8  & -1.01\%  & 166 \\
PAi-10 & 6.6  & 0.47\%  & 10.3 & -1.12\%  & 134 \\ \hline
\end{tabular}
\end{table}

\subsection{Solving the MIP model}\label{sec:mip-result}

In this section we investigate the performance of the mixed-integer programming algorithms over the proposed formulations. In particular, we examine the CPLEX MIP solver with default settings, with cutting plane algorithm disabled, and using automatic Benders decomposition algorithm. The tightness and computing time of a formulation's LP relaxation impact the performance of the B\&B algorithm. Additionally, the CPLEX heuristics' and the cutting plane algorithm's performance and the strength of  Benders cuts vary over the proposed formulations. We consider the FA and the DA formulations as existing models in the literature. Among the partial formulations, we employ PAe-5 and PAi-5 as they possess a satisfactory trade-off between the LP bound and computing time.   

The computing time for each instance is limited to 1 hour. In Tables~\ref{tab:mip-def}-\ref{tab:mip-benders-cut-long}, instances are divided into three groups: \emph{easy}, \emph{medium}, and \emph{difficult} instances. Easy instances are instances that all formulations have been able to prove optimality for within the time limit. Medium instances are instances for which some but not all formulations have proved optimality for. The other instances are classified as difficult. Note that the division of instances into the three groups and identifying the best solutions are based on the comparisons among the formulations under the corresponding settings of each table separately. Therefore, the difficulty classification of instances and the best solutions vary across tables. All numbers are reported as averages across instances in each group. The numbers in the parenthesis show the number of instances for each group. The notation used in the tables are as follows: \emph{Time}: solution time, \emph{\#n}: the number of explored nodes in the B\&B tree, \emph{\#B}: the number of instances for which the best solution is obtained by the corresponding formulation, \emph{\#O}: the number of instances where the corresponding formulation has been able to prove optimality within the time limit, \emph{B\%}: gap from the best solution, \emph{G\%}: optimality gap for the instances where their optimality has not been proven in the time limit. 

Table~\ref{tab:mip-def} reports the average performance of the CPLEX default MIP solver over the formulations per each instance class and overall. According to these results, the FA formulation outperforms other formulations over all instance groups and classes, in terms of proving the optimality and finding good-quality, feasible solutions. Analysis of the CPLEX logs shows that the heuristics perform better over the FA formulation. Moreover, although the FA formulation suffers from a weak LP bound, the cutting plane algorithm resolves this issue. These results are very surprising given that the latest study by \cite{Chouman2017} showed that the DA formulation outperforms the FA formulation. The DA formulation proves the optimality of easy instances in slightly shorter times than the partial formulations. However, the partial formulations perform better over the difficult instances, particularly to obtain good-quality solutions in comparison with the DA formulation. Moreover, this better performance is intenser for the large instances with many commodities, the C class.  As the trend for long instances is similar, the results are not reported separately. 

\begin{table}[htb!]
\caption{Performance of the MIP solver with default setting over different formulations for all instances}
\label{tab:mip-def}
\footnotesize
\begin{tabu}to \textwidth {X[3.2,c]|X[0.7,c]X[1.1,c]X[1.1,c]|X[0.7,c]X[1.1,c]X[1.1,c]|X[0.7,c]X[1.1,c]X[1.1,c]|X[0.7,c]X[1.1,c]X[1.1,c]} 
\hline
Instance   type              & \multicolumn{3}{c|}{DA}                                        & \multicolumn{3}{c|}{FA}                                                     & \multicolumn{3}{c|}{PAe-5}                               & \multicolumn{3}{c}{PAi-5}                            \\ \hline
\multicolumn{1}{l|}{}        & \multicolumn{2}{c}{Time}           & \#n                      & \multicolumn{2}{c}{Time}                        & \#n                      & \multicolumn{2}{c}{Time}              & \#n             & \multicolumn{2}{c}{Time}           & \#n             \\ \cline{2-13} 
Easy - C (12)                & \multicolumn{2}{c}{423.8}          & 5,043                    & \multicolumn{2}{c}{\textbf{310.8}}              & 15,449                   & \multicolumn{2}{c}{364.2}             & 5977            & \multicolumn{2}{c}{545.9}          & 11518           \\
Easy - C+ (8)                & \multicolumn{2}{c}{219.6}          & 1,846                    & \multicolumn{2}{c}{\textbf{207.7}}              & 1,881                    & \multicolumn{2}{c}{298.4}             & 2,151           & \multicolumn{2}{c}{223.9}          & 1,759           \\
Easy - R (131)               & \multicolumn{2}{c}{119.3}          & 1,843                    & \multicolumn{2}{c}{\textbf{87.0}}               & 4,891                    & \multicolumn{2}{c}{149.2}             & 3,443           & \multicolumn{2}{c}{151.1}          & 3,703           \\
\textit{Easy - all (151)}    & \multicolumn{2}{c}{\textit{148.8}} & \textit{2,097}           & \multicolumn{2}{c}{\textit{\textbf{111.1}}}     & \textit{5,570}           & \multicolumn{2}{c}{\textit{174.2}}    & \textit{3,576}  & \multicolumn{2}{c}{\textit{186.3}} & \textit{4,221}  \\ \hline
                             & \#B          & \#O                  & G\%                      & \#B                  & \#O                      & G\%                      & \#B                 & \#O             & G\%             & \#B            & \#O               & G\%             \\ \cline{2-13} 
Medium - C (3)               & \textbf{3}   & \textbf{1}           & \textbf{0.22\%}          & \textbf{3}           & \textbf{1}               & 0.47\%                   & \textbf{3}          & \textbf{1}      & 0.54\%          & 2              & 0                 & 0.64\%          \\
Medium - R (5)               & 4            & \textbf{3}           & \textbf{0.40\%}          & 4                    & \textbf{3}               & 1.23\%                   & \textbf{5}          & 1               & 0.71\%          & 4              & 1                 & 0.93\%          \\
\textit{Medium - all (8)}    & \textit{7}   & \textit{\textbf{4}}  & \textit{\textbf{0.31\%}} & \textit{7}           & \textit{\textbf{4}}      & \textit{0.84\%}          & \textit{\textbf{8}} & \textit{2}      & \textit{0.65\%} & \textit{6}     & \textit{1}        & \textit{0.81\%} \\ \hline
                             & \#B          & B\%                  & G\%                      & \#B                  & B\%                      & G\%                      & \#B                 & B\%             & G\%             & \#B            & B\%               & G\%             \\ \cline{2-13} 
Difficult - C (16)           & 3            & 2.40\%               & 3.54\%                   & \textbf{11}          & \textbf{0.02\%}          & \textbf{1.40\%}          & 2                   & 0.21\%          & 1.77\%          & 7              & 0.10\%            & 1.87\%          \\
Difficult - C+ (1)           & 0            & 1.59\%               & 5.77\%                   & 0                    & 0.57\%                   & 4.57\%                   & 0                   & 1.24\%          & 5.30\%          & \textbf{1}     & \textbf{0.00\%}   & \textbf{4.10\%} \\
Difficult - R (17)           & 6            & 0.21\%               & 2.20\%                   & \textbf{9}           & \textbf{0.06\%}          & \textbf{1.95\%}          & 3                   & 0.17\%          & 2.32\%          & 3              & 0.14\%            & 2.40\%          \\
\textit{Difficult- all (34)} & \textit{9}   & \textit{1.28\%}      & \textit{2.94\%}          & \textit{\textbf{20}} & \textit{\textbf{0.06\%}} & \textit{\textbf{1.76\%}} & \textit{5}          & \textit{0.22\%} & \textit{2.15\%} & \textit{11}    & \textit{0.12\%}   & \textit{2.20\%} \\ \hline
\end{tabu}
\end{table}

As Table~\ref{tab:mip-def} shows, the cutting plane algorithm plays a significant role in the performance of the MIP solver over a formulation. Therefore, to investigate the impact of the trade-off between the LP bound and computing time on the B\&B tree,  we disable the CPLEX cutting plane algorithm. Based on Table~\ref{tab:mip-no-cut}, the partial formulations outperforms the DA and FA formulations in general. PAi-5 performs better over easy instances, while PAe-5 provides smaller optimality gaps over medium and difficult instances. Although heuristics are still able to provide good-quality solutions for the FA formulation, the optimality gaps are significantly poorer due to its weak LP bound. Table~\ref{tab:mip-no-cut-long} presents the results over long instances. The superiority of the partial formulations is more significant over the large instances with many commodities (C class) and also long instances.

\begin{table}[htb!]
\caption{Performance of MIP solver with cutting plane algorithm disabled over different formulations for all instances}
\label{tab:mip-no-cut}
\footnotesize
\begin{tabu}to \textwidth {X[3.5,c]|X[0.7,c]X[1.1,c]X[1.1,c]|X[0.7,c]X[1.1,c]X[1.55,c]|X[0.7,c]X[1.1,c]X[1.1,c]|X[0.7,c]X[1.1,c]X[1.1,c]} 
\hline
Instance   type            & \multicolumn{3}{c|}{DA}                                       & \multicolumn{3}{c|}{FA}                                & \multicolumn{3}{c|}{PAe-5}                                        & \multicolumn{3}{c}{PAi-5}                                         \\ \hline
\multicolumn{1}{l|}{}      & \multicolumn{2}{c}{Time}          & \#n                      & \multicolumn{2}{c}{Time}           & \#n              & \multicolumn{2}{c}{Time}               & \#n                      & \multicolumn{2}{c}{Time}                        & \#n             \\ \cline{2-13} 
Easy - C (7)               & \multicolumn{2}{c}{1.2}           & 1,115                    & \multicolumn{2}{c}{341.1}          & 1,017,368        & \multicolumn{2}{c}{\textbf{1.1}}       & 989                      & \multicolumn{2}{c}{1.4}                         & 1564            \\
Easy - C+ (6)              & \multicolumn{2}{c}{\textbf{1.9}}  & 2,176                    & \multicolumn{2}{c}{2.2}            & 4,309            & \multicolumn{2}{c}{2.8}                & 2,403                    & \multicolumn{2}{c}{2.3}                         & 2,402           \\
Easy - R (102)             & \multicolumn{2}{c}{41.7}          & 2,326                    & \multicolumn{2}{c}{98.3}           & 242,590          & \multicolumn{2}{c}{35.0}               & 2,718                    & \multicolumn{2}{c}{\textbf{21.4}}               & 3,918           \\
\textit{Easy - all (115)}  & \multicolumn{2}{c}{\textit{37.2}} & \textit{2,244}           & \multicolumn{2}{c}{\textit{108.0}} & \textit{277,318} & \multicolumn{2}{c}{\textit{31.3}}      & \textit{2,596}           & \multicolumn{2}{c}{\textit{\textbf{19.1}}}      & \textit{3,695}  \\ \hline
                           & \#B         & \#O                  & G\%                      & \#B            & \#O                & G\%              & \#B                  & \#O             & G\%                      & \#B                  & \#O                      & G\%             \\ \cline{2-13} 
Medium - C (4)             & \textbf{4}  & \textbf{4}           & \textbf{--}              & \textbf{4}     & 0                  & 10.22\%          & \textbf{4}           & \textbf{4}      & \textbf{--}              & \textbf{4}           & 2                        & 0.83\%          \\
Medium - R (33)            & 32          & 32                   & 0.60\%                   & 31             & 0                  & 7.53\%           & \textbf{33}          & 27              & 1.23\%                   & \textbf{33}          & 24                       & 2.53\%          \\
\textit{Medium - all (37)} & \textit{36} & \textit{\textbf{36}} & \textit{\textbf{0.60\%}} & \textit{35}    & \textit{0}         & \textit{7.83\%}  & \textit{\textbf{37}} & \textit{31}     & \textit{1.23\%}          & \textit{\textbf{37}} & \textit{26}              & \textit{2.22\%} \\ \hline
                           & \#B         & B\%                  & G\%                      & \#B            & B\%                & G\%              & \#B                  & B\%             & G\%                      & \#B                  & B\%                      & G\%             \\ \cline{2-13} 
Difficult - C (20)              & 7           & 1.91\%               & 2.97\%                   & 8              & 0.14\%             & 12.56\%          & 9                    & 0.09\%          & \textbf{1.60\%}          & \textbf{11}          & \textbf{0.07\%}          & 2.35\%          \\
Difficult - C+ (3)              & 2           & 0.10\%               & \textbf{3.10\%}          & \textbf{3}     & \textbf{0.00\%}    & 3.96\%           & 2                    & 0.34\%          & 3.96\%                   & 2                    & 0.29\%                   & 3.87\%          \\
Difficult - R (18)              & 6           & 0.17\%               & \textbf{2.28\%}          & 7              & 0.11\%             & 10.86\%          & 6                    & 0.12\%          & 2.83\%                   & \textbf{10}          & \textbf{0.07\%}          & 3.53\%          \\
\textit{Difficult - all (41)}   & \textit{15} & \textit{1.02\%}      & \textit{2.67\%}          & \textit{18}    & \textit{0.12\%}    & \textit{11.18\%} & \textit{17}          & \textit{0.13\%} & \textit{\textbf{2.32\%}} & \textit{\textbf{23}} & \textit{\textbf{0.08\%}} & \textit{2.98\%} \\ \hline
\end{tabu}
\end{table}

\begin{table}[htb!]
\caption{Performance of MIP solver with cutting plane algorithm disabled over different formulations for long instances}
\label{tab:mip-no-cut-long}
\footnotesize
\begin{tabu}to \textwidth {X[3.5,c]|X[0.7,c]X[1.1,c]X[1.1,c]|X[0.7,c]X[1.1,c]X[1.1,c]|X[0.7,c]X[1.1,c]X[1.1,c]|X[0.7,c]X[1.1,c]X[1.1,c]} 
\hline
Instance   type          & \multicolumn{3}{c|}{DA}                        & \multicolumn{3}{c|}{FA}                         & \multicolumn{3}{c|}{PAe-5}                              & \multicolumn{3}{c}{PAi-5}                                        \\ \hline
                         & \#B        & \#O             & G\%             & \#B        & \#O             & G\%              & \#B        & \#O             & G\%                      & \#B                 & \#O                      & G\%             \\ \cline{2-13} 
Medium - R (2)           & 1          & \textbf{1}      & \textbf{0.60\%} & \textbf{2} & 0               & 10.46\%          & \textbf{2} & 0               & 1.09\%                   & \textbf{2}          & \textbf{1}               & 4.59\%          \\ \hline
                         & \#B        & B\%             & G\%             & \#B        & B\%             & G\%              & \#B        & B\%             & G\%                      & \#B                 & B\%                      & G\%             \\ \cline{2-13} 
Difficult - C (11)            & 3          & 3.35\%          & 4.27\%          & 3          & 0.19\%          & 13.20\%          & \textbf{5} & 0.12\%          & \textbf{1.73\%}          & 4                   & \textbf{0.05\%}          & 2.46\%          \\
Difficult - R (5)             & 0          & 0.30\%          & \textbf{2.66\%} & 1          & 0.15\%          & 12.14\%          & 1          & 0.12\%          & 2.82\%                   & \textbf{3}          & \textbf{0.11\%}          & 3.76\%          \\
\textit{Difficult - all (16)} & \textit{3} & \textit{2.40\%} & \textit{3.77\%} & \textit{4} & \textit{0.18\%} & \textit{12.87\%} & \textit{6} & \textit{0.12\%} & \textit{\textbf{2.07\%}} & \textit{\textbf{7}} & \textit{\textbf{0.07\%}} & \textit{2.87\%} \\ \hline
\end{tabu}
\end{table}

Benders decomposition algorithms have been successfully applied to MCND and its extensions as an effective solution algorithm \citep{COSTA20051429, costa_cordeau_gendron_2007, ZETINA2019311}.  In the next set of experiments, we evaluate the performance of the Benders decomposition algorithm over the proposed formulations. We use the CPLEX automatic Benders decomposition algorithm. Our aim is to provide indicative results on the effects of different level of aggregations on the performance of the Benders algorithm. To solve the problem of course specialized Benders algorithms would be more efficient. All settings are left at their default, except the presolve which is disabled. According to our experiments,  the presolve  destroys the special structure of the formulations and increases the LP solution time. Moreover, we add single-node cut-set constraints for the source and sink nodes to assist the feasibility of the master problem as it starts with only the constraint set~\eqref{eq:y_binary}. The added constraints are as~\eqref{eq:benders-ineq-1} and~\eqref{eq:benders-ineq-2}. These constraints improve the performance of the Benders algorithm over all formulations on average.

\begin{subequations}
\begin{align}
	& \sum_{j \in \mathcal{N}_i^+} u_{ij} y_{ij} \geq \sum_{k \in \mathcal{K}:o^k=i } d^k   &&\forall  i \in \mathcal{N} \label{eq:benders-ineq-1}\\
	& \sum_{j \in \mathcal{N}_i^-} u_{ji} y_{ji} \geq \sum_{k \in \mathcal{K}:s^k=i } d^k  &&\forall  i \in \mathcal{N} \label{eq:benders-ineq-2} 
\end{align}
\end{subequations}
Table~\ref{tab:mip-benders-cut} displays the summary results for the performance of the Benders algorithm over all formulations. The DA formulation outperforms other formulations on average. However, the PAe-5 formulation has still the best performance on the difficult instances of C class in terms of the quality of the solutions and the optimality gaps. The trend differs over the long instances as Table~\ref{tab:mip-benders-cut-long}. Over these instances, the PAe-5 performs the best in terms of all of the criteria. Moreover, the PAi-5 outperforms the DA and FA formulations. This particularly emphasizes the effectiveness of the formulations based on partial aggregations as their impact is more significant over long instances to reduce the LP computing time. The results indicate that the subproblems of the partial formulations benefits from shorter computing times, yet generate strong Benders cuts. 

\begin{table}[htb!]
\caption{Performance of the Benders algorithm over different formulations with added single-node cut-set constraints~\eqref{eq:benders-ineq-1} and~\eqref{eq:benders-ineq-2} for all instances}
\label{tab:mip-benders-cut}
\footnotesize
\begin{tabu}to \textwidth {X[3.7,c]|X[0.7,c]X[1.1,c]X[1.4,c]|X[0.7,c]X[1.1,c]X[1.4,c]|X[0.7,c]X[1.1,c]X[1.4,c]|X[0.7,c]X[1.1,c]X[1.4,c]} 
\hline
Instance   type            & \multicolumn{3}{c|}{DA}                                                & \multicolumn{3}{c|}{FA}                               & \multicolumn{3}{c|}{PAe-5}                               & \multicolumn{3}{c}{PAi-5}                           \\ \hline
\multicolumn{1}{l}{}       & \multicolumn{2}{|c}{Time}                    & \#n                      & \multicolumn{2}{c}{Time}           & \#n              & \multicolumn{2}{c}{Time}               & \#n             & \multicolumn{2}{c}{Time}          & \#n             \\ \cline{2-13} 
Easy - C+ (5)              & \multicolumn{2}{c}{\textbf{137.7}}          & 210,528                  & \multicolumn{2}{c}{540.4}          & 539,608          & \multicolumn{2}{c}{191.3}              & 255,806         & \multicolumn{2}{c}{433.6}         & 375,894         \\
Easy - R (75)              & \multicolumn{2}{c}{\textbf{28.5}}           & 11,424                   & \multicolumn{2}{c}{130.3}          & 115,893          & \multicolumn{2}{c}{39.6}               & 21,475          & \multicolumn{2}{c}{39.3}          & 17,111          \\
\textit{Easy - all (80)}   & \multicolumn{2}{c}{\textit{\textbf{35.3}}}  & \textit{23,868}          & \multicolumn{2}{c}{\textit{156.0}} & \textit{142,375} & \multicolumn{2}{c}{\textit{49.1}}      & \textit{36,121} & \multicolumn{2}{c}{\textit{63.9}} & \textit{39,535} \\ \hline
                           & \#B                  & \#O                  & G\%                      & \#B           & \#O                & G\%              & \#B         & \#O                      & G\%             & \#B           & \#O               & G\%             \\ \cline{2-13} 
Medium - C (8)             & \textbf{8}           & \textbf{8}           & \textbf{--}              & 1             & 0                  & 5.42\%           & 7           & 7                        & 0.27\%          & 7             & 7                 & 0.57\%          \\
Medium - C+ (1)            & \textbf{1}           & \textbf{1}           & \textbf{--}              & 1             & 0                  & 0.43\%           & \textbf{1}  & \textbf{1}               & \textbf{--}     & \textbf{1}    & \textbf{1}        & \textbf{--}     \\
Medium - R (29)            & \textbf{29}          & \textbf{29}          & \textbf{--}              & 0             & 0                  & 9.52\%           & 26          & 24                       & 1.66\%          & 24            & 23                & 1.52\%          \\
\textit{Medium - all (36)} & \textit{\textbf{38}} & \textit{\textbf{38}} & \textit{\textbf{--}}     & \textit{2}    & \textit{0}         & \textit{8.42\%}  & \textit{34} & \textit{32}              & \textit{1.43\%} & \textit{32}   & \textit{31}       & \textit{1.39\%} \\ \hline
                           & \#B                  & B\%                  & O\%                      & \#B           & B\%                & O\%              & \#B         & B\%                      & O\%             & \#B           & B\%               & O\%             \\ \cline{2-13}
Difficult - C (23)              & \textbf{13}          & 3.67\%               & 6.54\%                   & 0             & 5.69\%             & 23.59\%          & 4           & \textbf{0.73\%}          & \textbf{5.43\%} & 6             & 1.03\%            & 6.33\%          \\
Difficult - C+ (3)              & \textbf{1}           & \textbf{0.14\%}      & \textbf{11.13\%}         & \textbf{1}    & 6.89\%             & 16.74\%          & 0           & 5.06\%                   & 14.52\%         & \textbf{1}    & 6.68\%            & 15.71\%         \\
Difficult - R (49)              & \textbf{28}          & \textbf{0.63\%}      & \textbf{5.61\%}          & 0             & 4.70\%             & 21.69\%          & 16          & 0.90\%                   & 7.36\%          & 7             & 1.32\%            & 8.19\%          \\
\textit{Difficult - all (75)}   & \textit{\textbf{42}} & \textit{1.54\%}      & \textit{\textbf{6.11\%}} & \textit{1}    & \textit{5.09\%}    & \textit{22.08\%} & \textit{20} & \textit{\textbf{1.01\%}} & \textit{7.06\%} & \textit{14}   & \textit{1.45\%}   & \textit{7.92\%} \\ \hline
\end{tabu}
\end{table}

\begin{table}[htb!]
\caption{Performance of the Benders algorithm over different formulations with added single-node cut-set constraints~\eqref{eq:benders-ineq-1} and~\eqref{eq:benders-ineq-2} for long instances}
\label{tab:mip-benders-cut-long}
\footnotesize
\begin{tabu}to \textwidth {X[3.7,c]|X[0.7,c]X[1.1,c]X[1.4,c]|X[0.7,c]X[1.1,c]X[1.4,c]|X[0.7,c]X[1.1,c]X[1.4,c]|X[0.7,c]X[1.1,c]X[1.4,c]} 
\hline
Instance   type          & \multicolumn{3}{c|}{DA}                         & \multicolumn{3}{c|}{FA}                         & \multicolumn{3}{c|}{PAe-5}                                                & \multicolumn{3}{c}{PAi-5}                     \\ \hline
                         & \#B        & B\%             & O\%             & \#B        & B\%             & O\%              & \#B                 & B\%                      & O\%                      & \#B        & B\%             & O\%             \\ \cline{2-13} 
Difficult - C (11)            & 4          & 7.66\%          & 9.84\%          & 0          & 6.05\%          & 23.81\%          & 2                   & \textbf{0.51\%}          & \textbf{5.67\%}          & \textbf{5} & 0.60\%          & 6.28\%          \\
Difficult - R (7)             & 0          & 2.58\%          & 8.89\%          & 0          & 6.29\%          & 27.89\%          & \textbf{5}          & \textbf{0.58\%}          & \textbf{8.38\%}          & 2          & 2.53\%          & 10.93\%         \\
\textit{Difficult - all (18)} & \textit{4} & \textit{5.68\%} & \textit{9.47\%} & \textit{0} & \textit{6.14\%} & \textit{25.40\%} & \textit{\textbf{7}} & \textit{\textbf{0.54\%}} & \textit{\textbf{6.72\%}} & \textit{7} & \textit{1.35\%} & \textit{8.09\%} \\ \hline
\end{tabu}
\end{table}

\begin{figure}[htb!]
\begin{subfigure}{.5\textwidth}
  \centering
  \includegraphics[width=1\linewidth]{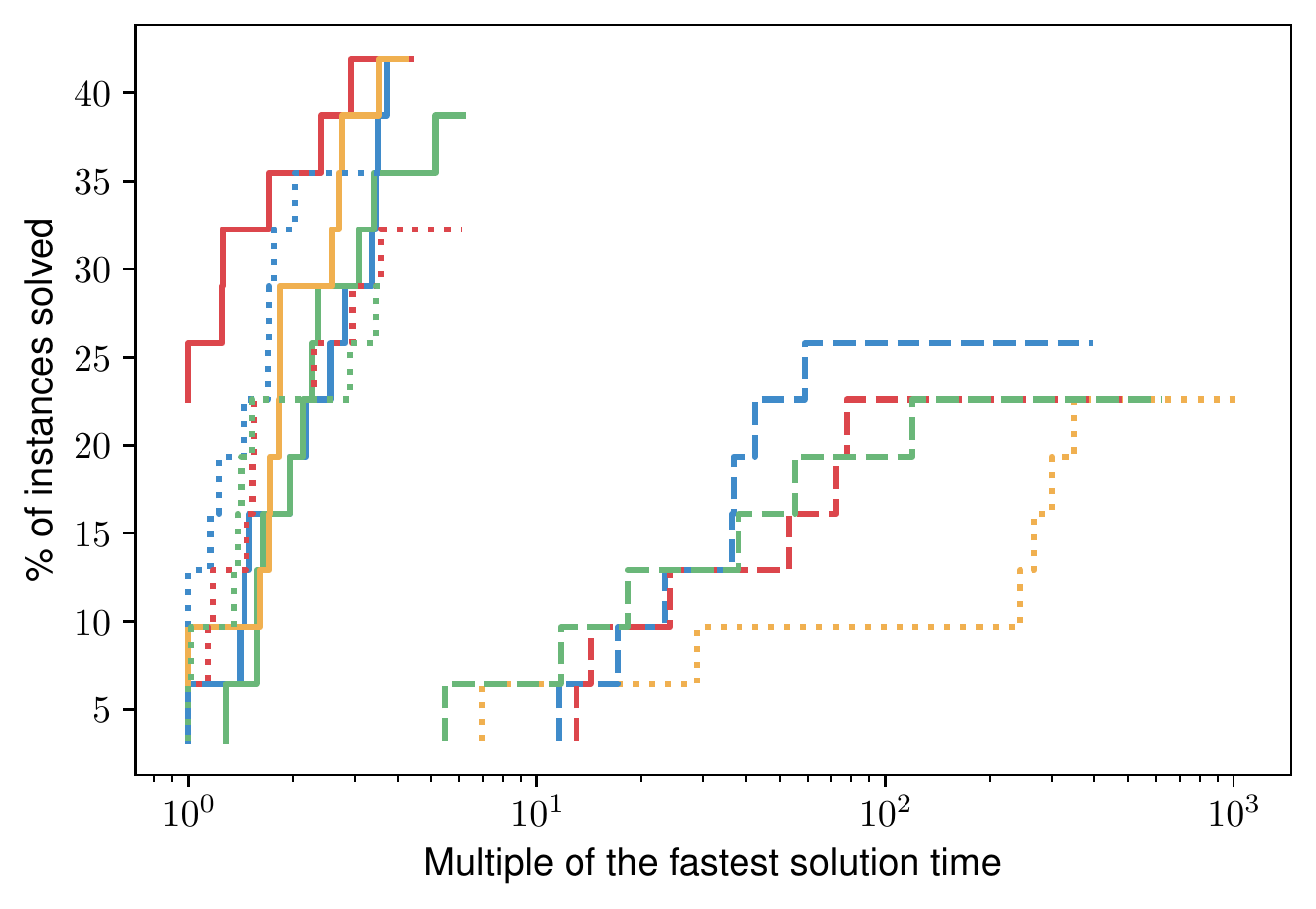}  
  \caption{Class C}
  \label{fig:profile_C}
\end{subfigure}
\begin{subfigure}{.5\textwidth}
  \centering
  \includegraphics[width=1\linewidth]{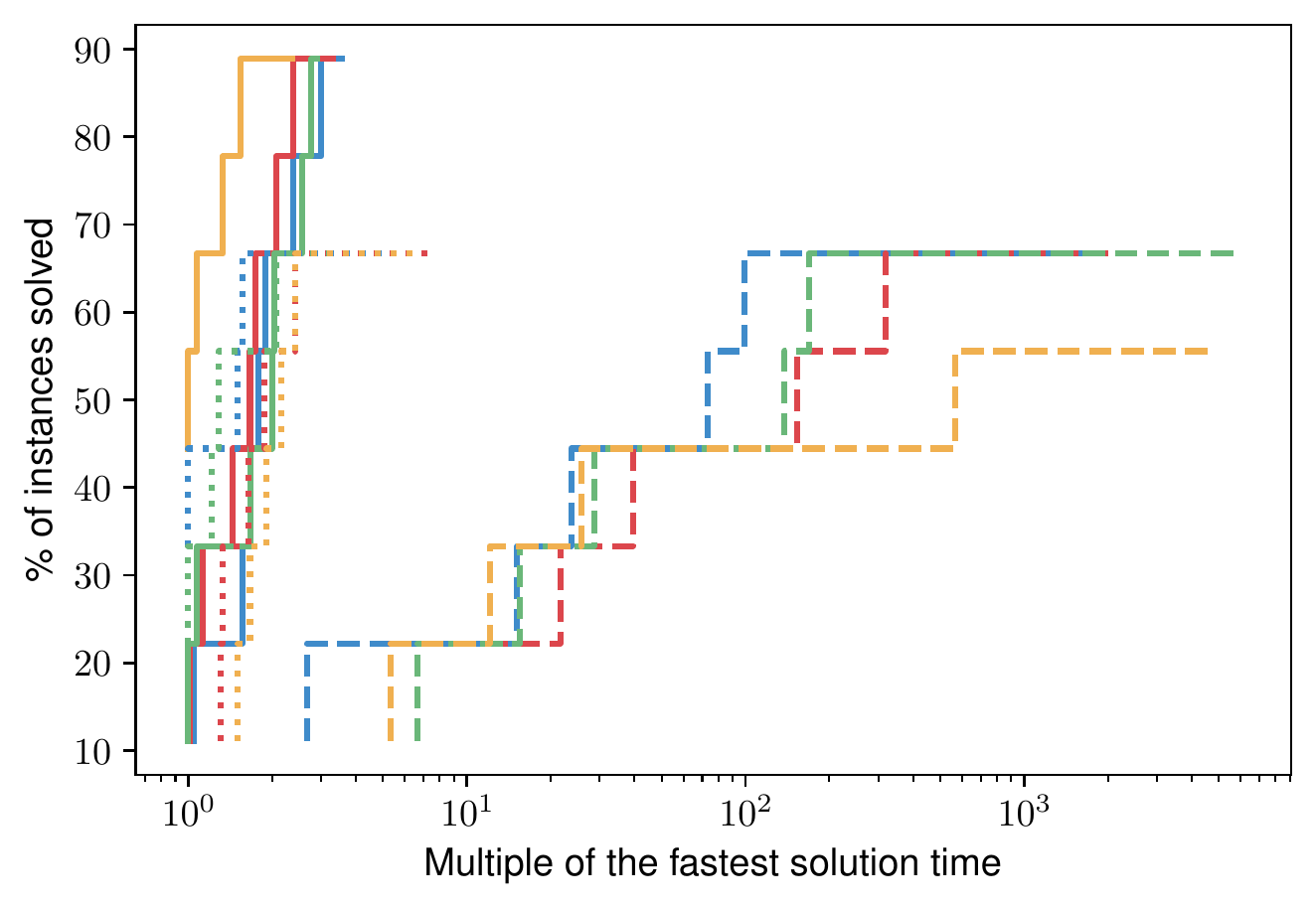}  
  \caption{Class C+}
  \label{fig:profile_C+}
\end{subfigure}


\begin{subfigure}{.5\textwidth}
  \centering
  \includegraphics[width=1\linewidth]{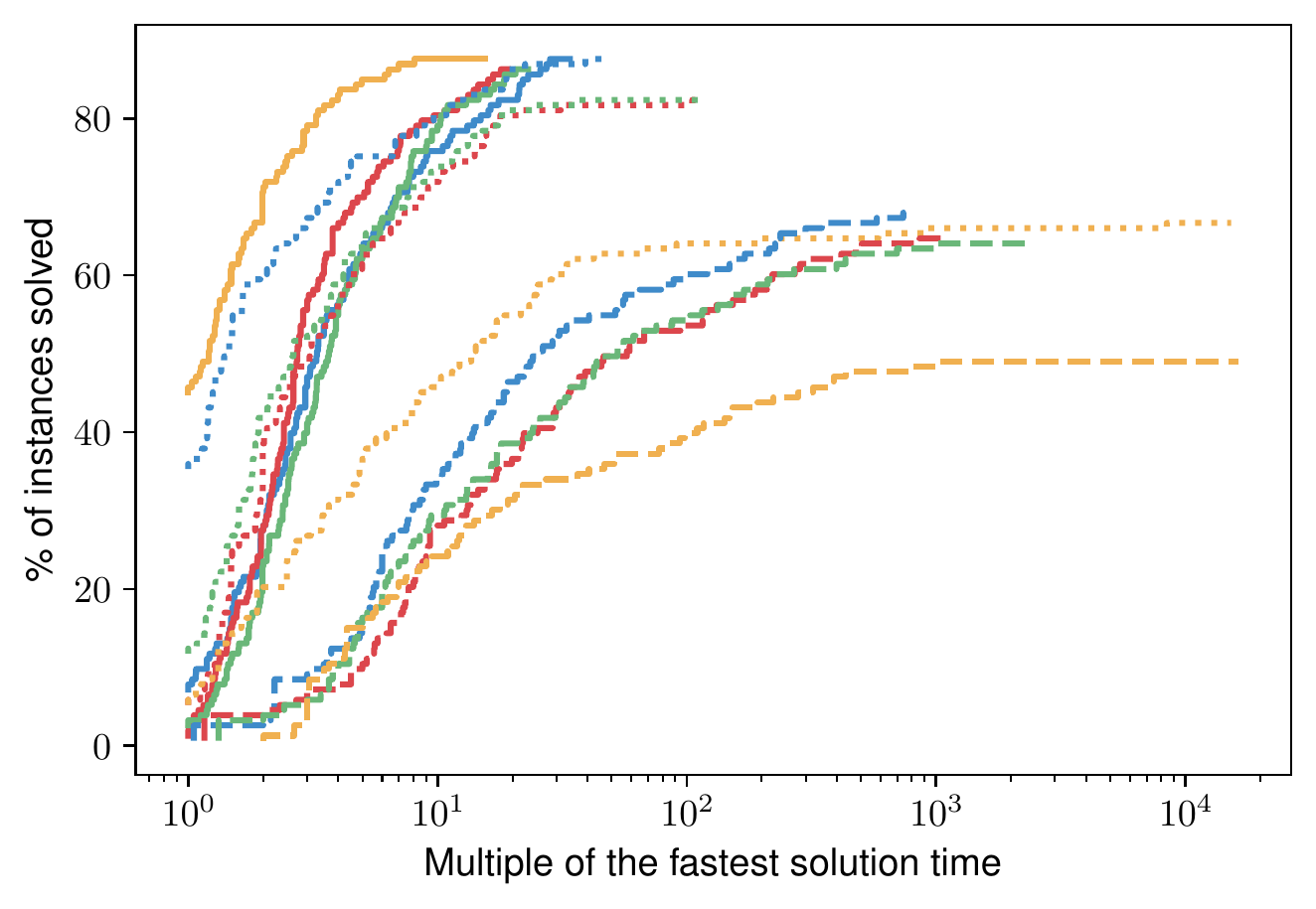}  
  \caption{Class R}
  \label{fig:profile_R}
\end{subfigure}
\begin{subfigure}{.5\textwidth}
  \centering
  \includegraphics[width=1\linewidth]{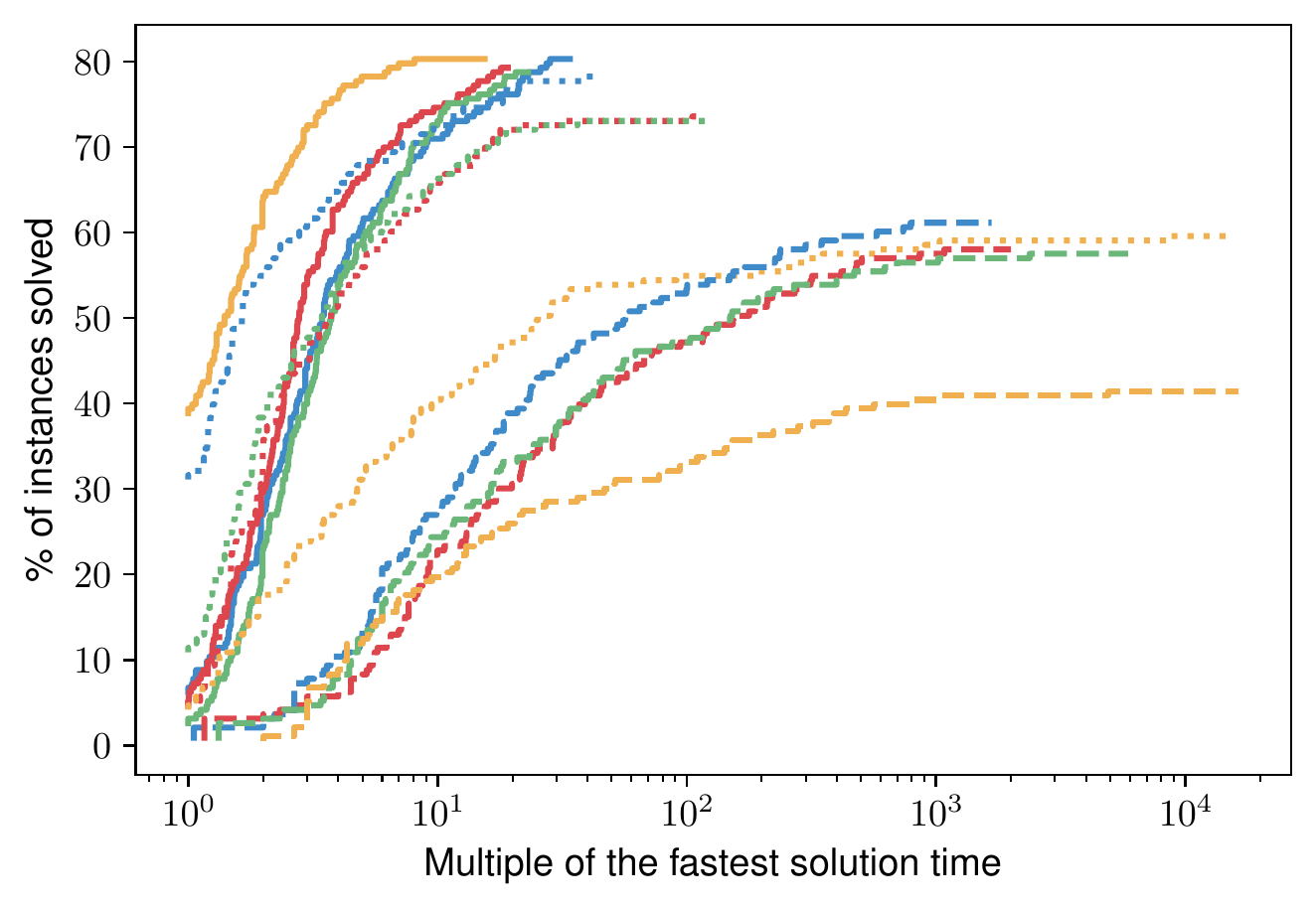}  
  \caption{All instances}
  \label{fig:profile_all}
\end{subfigure}
\begin{subfigure}{1\textwidth}
  \centering
  \includegraphics[width=0.5\linewidth]{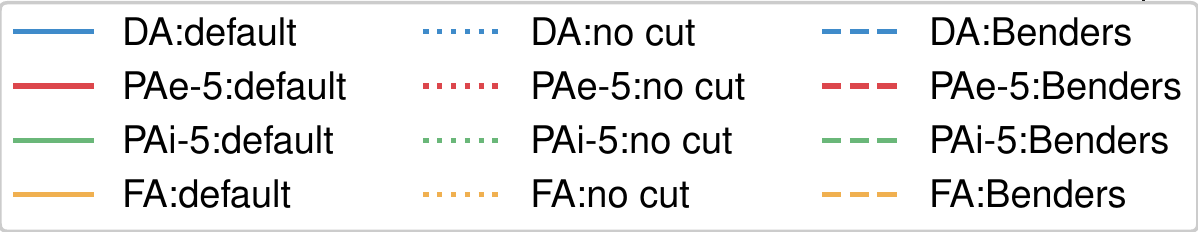}
\end{subfigure}

\caption{Performance profile of the formulations based on different aggregation schemes with respect to solution time (log scale)}
\label{fig:profiles_compare}
\end{figure}

The summary of the results on solving the MIP model is presented as performance profiles in Figure~\ref{fig:profiles_compare}. Such performance profiles are proposed by \cite{dolan_more_2002} to evaluate and compare optimization solvers. A performance profile shows the percentage of instances solved to the optimality within the multiple time of the fastest solution time of each instance. Performance profiles are reported per instance class and overall. Generally, the default setting performs better as expected. For the class C of instances, the PAe-5 formulation performs the best. This formulation proves the most percentage of these instances in less time in comparison with other formulations. Moreover, it solves more instances in the fastest time in comparison with other formulations. In other words, PAe-5 using CPLEX with the default settings is able to prove the optimality of almost 50\% of instances, which their optimality is proven in the time limit, in the fastest time. This is particularly interesting since this class of instances include large instances with many commodities, and partial aggregation is effective on such instances.  The trend for other instance classes and average indicate the superiority of the FA formulation by using CPLEX with the default settings as discussed before. Benders decomposition and FA with cutting plane algorithm disabled performs poorly in comparison with other formulations. However, our goal for Benders decomposition experiments is to study the impacts of aggregation schemes on its performance and not to  compare against state-of-the-art technology. While the cutting plane algorithm significantly improves the performance of the FA formulation, the DA formulation with cutting plane algorithm disabled performs better than the DA with default settings. In general, in shorter computing times, the DA with cutting plane algorithm disabled perform better than the PAe-5 with the default setting. However, the performance profile of the PAe-5 with the default setting surpasses the DA with cutting plane algorithm disabled as the time passes.

\section{Conclusions and future research}\label{sec:conclusion}

This paper introduces new commodity representations for multicommodity network flow problems. Such commodity representations are utilized to aggregate constraints and variables in the MIP model partially instead of the conventional full aggregation. We apply the proposed partial aggregations to the multicommodity fixed-charge network design problem to show the implications of this approach. However, they are applicable to any other problem that can be formulated with a multicommodity network flow subproblem, and for which there exist an alternative fully-aggregated formulation. We propose two variants of MIP formulations based on this partial aggregation concept and a heuristic algorithm for constructing partial aggregations of various sizes. These MIPs are proven to provide a valid lower bound MCND which is tighter than the fully aggregated formulation. The proposed MIPs and their LP relaxations are evaluated empirically through an extensive computational study on benchmark instances. The results indicate that the partial formulations have significantly shorter computing times in comparison with the DA formulation  with only small loss in the LP bound in contrast to the FA formulation. The proposed formulations have dimensions in between the dimensions of the DA and FA formulations. Furthermore, the proposed models form a Pareto frontier for the trade-off between the LP bound and computing time, whether only the LP relaxation is solved or a cutting plane algorithm is employed. We investigate the performance of the mixed-integer programming algorithms over the formulations.  The PAe-5 formulation is particularly effective for the large instances with many commodities and long instances. In addition, the proper trade-off between the LP bound and computing time benefits the partial formulations in the B\&B algorithm. The partial formulations are also beneficial for long instances using a Benders decomposition algorithm since the reduced sub-problems have shorter computing times yet generate effective cuts.

We propose three possible directions for future research. First, partially-aggregated formulations can be employed to develop specialized solution algorithms. Particularly, our early experiments show that the Benders decomposition algorithm performs well over the partial formulations for the long instances and large instances with many commodities. The second direction is to develop valid cuts for the partial formulations. The DA and FA formulations have been extensively studied to develop efficient cutting plane algorithms. Similarly, cutting plane algorithms would be beneficial to solvers using partial aggregations. The third direction is to apply this partial aggregation approach to other multicommodity network flow problems and study the implications. The extensions are particularly straightforward for the problems where there is no commodity-specific attribute over the underlying network.

\section*{Acknowledgments}

This research is supported by Australian Research Council with Grant LP160100547. This research is also supported in part by the Monash eResearch Centre through the use of the MonARCH HPC Cluster.

\small
\bibliography{mybibfile}

\end{document}